\documentclass[final,onefignum,onetabnum]{siamart190516}

\usepackage{bm}
\usepackage{mathrsfs}
\usepackage{amsmath}
\usepackage{amsfonts}
\usepackage{mathtools}

\usepackage{epstopdf}

\usepackage[english]{babel}
\usepackage{cleveref}
\usepackage{graphicx}
\usepackage{subcaption}
\usepackage{wrapfig}
\usepackage{float}
\usepackage{color}
\usepackage{verbatim}
\usepackage{algorithmic}
\usepackage{multirow}
\usepackage[titletoc, title]{appendix}

\usepackage{wasysym}

\usepackage{amsopn}

\newsiamremark{remark}{Remark}
\newsiamremark{assumption}{Assumption}

\newcommand{\numberthis}{\addtocounter{equation}{1}\tag{\theequation}}
\newcommand{\del}{\partial}

\title{Adaptive Deep Learning for High-Dimensional Hamilton-Jacobi-Bellman Equations\thanks{Preliminary results of this work appeared in the American Control Conference 2020 \cite{Nakamura2020}. {\funding: The work of the first and second authors was partially supported with funding from the Defense Advanced Research Projects Agency (DARPA) grant FA8650-18-1-7842.}
}
}

\author{Tenavi Nakamura-Zimmerer\thanks{Department of Applied Mathematics, Baskin School of Engineering, University of California, Santa Cruz (\email{tenakamu@ucsc.edu}, \email{qgong@ucsc.edu}).}
\and
Qi Gong\footnotemark[2]
\and
Wei Kang\thanks{Department of Applied Mathematics, Naval Postgraduate School, Monterery, CA, (\email{wkang@nps.edu}).}}

\headers{Deep Learning for HJB Equations}{T. Nakamura-Zimmerer, Q. Gong, and W. Kang}

\ifpdf
  \DeclareGraphicsExtensions{.eps,.pdf,.png,.jpg}
  \hypersetup{
  pdftitle={Adaptive Deep Learning for High-Dimensional Hamilton-Jacobi-Bellman Equations},
  pdfauthor={T. Nakamura-Zimmerer, Q. Gong, and W. Kang}
}
\else
  \DeclareGraphicsExtensions{.eps}
\fi

\begin{document}

\maketitle

\begin{abstract}
Computing optimal feedback controls for nonlinear systems generally requires solving Hamilton-Jacobi-Bellman (HJB) equations, which are notoriously difficult when the state dimension is large. Existing strategies for high-dimensional problems often rely on specific, restrictive problem structures, or are valid only locally around some nominal trajectory. In this paper, we propose a data-driven method to approximate semi-global solutions to HJB equations for general high-dimensional nonlinear systems and compute candidate optimal feedback controls in real-time. To accomplish this, we model solutions to HJB equations with neural networks (NNs) trained on data generated without discretizing the state space. Training is made more effective and data-efficient by leveraging the known physics of the problem and using the partially-trained NN to aid in adaptive data generation. We demonstrate the effectiveness of our method by learning solutions to HJB equations corresponding to the attitude control of a six-dimensional nonlinear rigid body, and nonlinear systems of dimension up to 30 arising from the stabilization of a Burgers'-type partial differential equation. The trained NNs are then used for real-time feedback control of these systems.
\end{abstract}



\section{Introduction}

For the optimal control of nonlinear dynamical systems, it is well-known that open-loop controls are not robust to model uncertainty or disturbances. For slowly evolving processes, it is possible to use model predictive control by recomputing the open-loop optimal solutions for relatively short time horizons in the future. However, for most applications one typically desires a feedback control law, as feedback controls are inherently more robust to disturbances. In principle, optimal feedback controllers can synthesized by solving a (discretized) Hamilton-Jacobi-Bellman (HJB) equation, a partial differential equation (PDE) in $n$ spatial dimensions plus time. The size of the discretized problem increases exponentially with $n$, making direct solution intractable for even moderately large problems. This is the so-called ``curse of dimensionality.''

For this reason, there is an extensive literature on methods of finding approximate solutions for HJB equations. Some key examples include series expansions \cite{Albrekht1961, Lukes1969, Kang1992}, level set methods \cite{Osher1988}, patchy dynamic programming \cite{Cacace2012, Navasca2007}, semi-Lagrangian methods \cite{Bokanowski2013, Falcone2013}, method of characteristics and Hopf formula-based algorithms \cite{Darbon2016, Chow2019, Yegorov2018}, and polynomial approximation \cite{Kalise2018}. These existing methods suffer one or more of the following drawbacks: the problem's dimension is limited; the accuracy of the solution is hard to verify for general systems; the solution may be valid only in a small region; or the system model must have certain special algebraic structure.

In \cite{Kang2015, Kang2017_COA}, semi-global solutions to HJB equations are computed by combining the method of characteristics with sparse state space discretization. In this approach, a two-point boundary value problem (BVP) is solved at each point in a sparse grid. These BVPs can be solved independently, making the algorithm \textit{causality-free}. This property is attractive because the computation does not depend on a grid, and hence they can be applied to high-dimensional problems. The Hopf formula methods \cite{Darbon2016, Chow2019, Yegorov2018} also have this property, though it is achieved in a different way and under certain convexity/concavity assumptions. Causality-free methods are usually too slow for online computation, but they are perfectly parallelizable so can be used to generate large data sets offline. Such data sets can then be used to construct faster solutions such as sparse grid interpolants \cite{Kang2015, Kang2017_COA} or, as in this paper, neural networks.

Using neural networks (NNs) as a basis for solving HJB equations is not by itself a new idea, and deep learning approaches have led to promising results; see for instance \cite{Khalaf2005, Cheng2007, Tassa2007, Jiang2016, Sirignano2018, Izzo2019}. To the best of our knowledge, state-of-the-art NN-based techniques generally rely on either minimizing the residual of the PDE and (artificial) boundary conditions at randomly sampled collocation points \cite{Khalaf2005, Cheng2007, Tassa2007, Sirignano2018}; or, due to computational limitations, approximating the control and/or HJB solution and its gradient in a small neighborhood of a nominal trajectory \cite{Jiang2016, Izzo2019}. In \cite{Darbon2020, Darbon2021}, a specialized NN architecture is proposed to solve some classes of Hamilton-Jacobi equations, but this method has yet to be generalized to state-dependent HJB equations arising in optimal control. Deep learning techniques have also been proposed for solving high-dimensional stochastic optimal control problems (see e.g. \cite{Han2018_PNAS, Hure2018, Bachouch2018}).

In this paper, we develop a computational method for solving high-dimensional HJB equations and synthesizing candidate optimal feedback controllers. Our approach is data-driven and consists of three main steps. First, we generate a small set of open-loop optimal control solutions using a causality-free algorithm based on Pontryagin's Minimum Principle (PMP). In the second step, we use the data set to train a NN to approximate the solution to the HJB equation, called the \textit{value function}. During training, we supply information about the value function gradient, which encourages the NN to learn the shape of the value function rather than just fitting point data. We also estimate the number of samples needed to obtain a good model. Additional samples are chosen in regions where the value function is difficult to learn, and are obtained quickly with the aid of the NN. In this sense our method involves \textit{adaptive sampling}. Lastly, the accuracy of the NN is verified on independent data generated using the same causality free algorithm from the first step. Unlike other NN-based methods for deterministic HJB equations, our approach does not require computing expensive PDE residuals and the solution is valid over large spatial domains.

As an illustrative example, the method is applied to design an attitude controller of a rigid-body satellite equipped with momentum wheels. This is a highly nonlinear problem with $n=6$ spatial dimensions and $m=3$ control inputs. With the proposed method, we obtain a model of the value function with accuracy comparable to that obtained in \cite{Kang2017_COA}, but require far fewer sample trajectories to do so. Scalability of the method is tested on problems of dimension $n = 10$, 20, and 30 arising from pseudospectral discretization of a Burgers'-type PDE. We show that the method is capable of handling these high-dimensional problems without simplifying the dynamics.

Through these examples, we demonstrate several advantages and potential capabilities of the proposed framework. These include solving HJB equations over semi-global domains with empirically validated levels of accuracy, progressive generation of rich data sets, and computationally efficient nonlinear feedback control for real-time applications. Solution of high-dimensional problems is enabled by efficient and adaptive causality-free data generation, physics-informed learning, and the inherent capacity of NNs for dealing with high-dimensional data.

\subsection{Abbreviations and notation}

Here we present a brief list of some of the abbreviations, terminology, and notation used in this paper.

\vspace{10pt}
\begin{center}
\begin{tabular}{lcl}
\hline
OCP & \dots & optimal control problem\\
HJB & \dots & Hamilton-Jacobi-Bellman equation\\
PMP & \dots & Pontryagin's Minimum Principle\\
NN & \dots & neural network\\
RMAE & \dots & relative mean absolute error\\
RM$L^2$ & \dots & relative mean $L^2$ error\\
$\mathcal D$ & \dots & data set\\
$\mu$ & \dots & gradient regularization weight \\
$C$ & \dots & adaptive sampling convergence parameter \\
\hline
\end{tabular}
\end{center}
\vspace{10pt}

\section{A causality-free method for HJB equations}
\label{sec: BVP}

We consider fixed final time optimal control problems (OCP) of the form
\begin{equation}
\label{eq: OCP}
\left \{
\begin{array}{cl}
\underset{\bm u (\cdot) \in \mathbb U}{\text{minimize}} & J \left[ \bm u (\cdot) \right] = F (\bm x (t_f)) + \displaystyle \int_0^{t_f} \mathcal L (t, \bm x, \bm u) dt , \\
\text{subject to} & \dot {\bm x} (t) = \bm f (t, \bm x, \bm u) , \\
	& \bm x (0) = \bm x_0 .
\end{array}
\right .
\end{equation}
Here $\bm x (t) : [0, t_f] \to \mathbb X \subseteq \mathbb R^n$ is the state, $\bm u (t, \bm x) : [0, t_f] \times \mathbb X \to \mathbb U \subseteq \mathbb R^m$ is the control, and $\bm f (t, \bm x, \bm u) : [0, t_f] \times \mathbb X \times \mathbb U \to \mathbb R^n$ is a Lipschitz continuous vector field. $J \left[ \bm u (\cdot) \right]$ is the cost functional which is composed of $F (\bm x (t_f)) : \mathbb X \to \mathbb R$, the terminal cost, and $\mathcal L (t, \bm x, \bm u) : [0, t_f] \times \mathbb X \times \mathbb U \to \mathbb R$, the running cost. We assume that the cost functional is convex in $\bm x$ and $\bm u$. In this paper we consider the case where the final time $t_f < \infty$ is fixed.

For a given initial condition $\bm x (0) = \bm x_0$, many numerical methods exist to compute the optimal open-loop solution,
\begin{equation}
\label{eq: optimal open-loop}
\bm u = \bm u^* (t; \bm x_0) .
\end{equation}
The open-loop control \cref{eq: optimal open-loop} which solves \cref{eq: OCP} is valid for all $t \in [0, t_f]$, but only for the fixed initial condition $\bm x (0) = \bm x_0$. Due to various sources of disturbance and real-time application requirements, for practical implementation one typically desires an optimal control in closed-loop feedback form,
\begin{equation}
\label{eq: optimal feedback}
\bm u = \bm u^* (t, \bm x) ,
\end{equation}
which can be evaluated online given any $t \in [0, t_f]$ and a measurement of $\bm x \in \mathbb X$.

To compute the optimal feedback control, we follow the standard procedure in dynamic programming (see e.g. \cite{Liberzon2011}) and define the value function $V(t,\bm x) : [0, t_f] \times \mathbb X \to \mathbb R$ as the optimal cost-to-go of \cref{eq: OCP} starting at $(t, \bm x)$. That is,
\begin{equation}
\label{eq: value function}
V (t, \bm x) \coloneqq J \left[ \bm u^* (\cdot) \right] =
\left \{
\begin{array}{cl}
\displaystyle \inf_{\bm u (\cdot) \in \mathbb U} & F (\bm y (t_f)) + \displaystyle \int_t^{t_f} \mathcal L (\tau, \bm y, \bm u) d\tau , \\
\text{s.t.} & \dot {\bm y} (\tau) = \bm f (\tau, \bm y, \bm u) , \\
	& \bm y (t) = \bm x .
\end{array}
\right .
\end{equation}
It can be shown that the value function is the unique viscosity solution \cite{Crandall1983} of the Hamilton-Jacobi-Bellman (HJB) PDE,
\begin{equation}
\label{eq: HJB1}
\begin{dcases*}
- V_t (t,\bm x) - \min_{\bm u \in \mathbb U} \left\{ \mathcal L (t, \bm x, \bm u) + [V_{\bm x} (t,\bm x)]^T \bm f (t, \bm x, \bm u) \right\} = 0 , \\
V(t_f,\bm x) = F (\bm x) ,
\end{dcases*}
\end{equation}
where we denote $V_t \coloneqq \del V / \del t$ and $V_{\bm x} \coloneqq \left[ \del V / \del \bm x \right]^T$. Note that if the value function is $C^2$, then it is the unique classical solution of \cref{eq: HJB1}.

To compute the control given the value function $V (\cdot)$, we start by defining the Hamiltonian
\begin{equation}
\label{eq: Hamiltonian}
\mathcal H (t, \bm x, \bm \lambda, \bm u) \coloneqq\mathcal L (t, \bm x, \bm u) + \bm \lambda^T \bm f (t, \bm x, \bm u) ,
\end{equation}
where $\bm \lambda (t) : [0, t_f] \to \mathbb R^n$ is the costate. The optimal control satisfies the Hamiltonian minimization condition,
\begin{equation}
\label{eq: optimal control as Hamiltonian minimizer}
\bm u^* (t) = \bm u^* (t, \bm x; \bm \lambda) = \underset{\bm u \in \mathbb U}{\text{arg min}} \, \mathcal H (t, \bm x, \bm \lambda, \bm u) .
\end{equation}
If we denote the minimized Hamiltonian by $\mathcal H^* (t, \bm x, \bm \lambda) \coloneqq \mathcal H \left( t, \bm x, \bm \lambda, \bm u^* (t, \bm x; \bm \lambda) \right)$, then \cref{eq: HJB1} can be expressed as
\begin{equation}
\label{eq: HJB}
\begin{dcases*}
- V_t (t, \bm x) - \mathcal H^* \left( t, \bm x, V_{\bm x} \right) = 0 , \\
V (t_f, \bm x) = F (\bm x) .
\end{dcases*}
\end{equation}
If \cref{eq: HJB} can be solved (in the viscosity sense), then it provides both necessary and sufficient conditions for optimality. Moreover, the optimal feedback control is computed by substituting
\begin{equation}
\label{eq: costate as gradient}
\bm \lambda (t) = V_{\bm x} (t, \bm x)
\end{equation}
into \cref{eq: optimal control as Hamiltonian minimizer} to get
\begin{equation}
\label{eq: optimal feedback control using dVdx}
\bm u^* (t, \bm x) = \bm u^* (t, \bm x; V_{\bm x}) = \underset{\bm u \in \mathbb U}{\text{arg min}} \, \mathcal H \left( t, \bm x, V_{\bm x}, \bm u \right) .
\end{equation}
This means that with $V_{\bm x} (\cdot)$ available, the feedback control is obtained as the solution of an (ideally straightforward) optimization problem.

\subsection{Pontryagin's Minimum Principle}

To make use of \cref{eq: optimal feedback control using dVdx}, we need an efficient way to approximate the value function and its gradient. Like \cite{Kang2015, Kang2017_COA}, rather than solve the full HJB equation \cref{eq: HJB} on a grid, we exploit the fact that the characteristics of solutions to \cref{eq: HJB} evolve according to a two-point BVP, well-known in optimal control as Pontryagin's Minimum Principle (PMP):
\begin{equation}
\label{eq: BVP}
\begin{dcases*}
\dot{\bm x} (t) = \mathcal H_{\bm \lambda} = \bm f (t, \bm x, \bm u^* (t, \bm x; \bm \lambda)) ,
	& \qquad $\bm x (0) = \bm x_0$, \\
\dot{\bm \lambda} (t) = - \mathcal H_{\bm x} (t, \bm x, \bm \lambda, \bm u^* (t, \bm x; \bm \lambda)) ,
	& \qquad $\bm \lambda (t_f) = F_{\bm x} (\bm x (t_f))$, \\
\dot v (t) = - \mathcal L (t, \bm x, \bm u^* (t, \bm x; \bm \lambda)) ,
	& \qquad $v (t_f) = F (\bm x (t_f))$.
\end{dcases*}
\end{equation}
The two-point BVP provides a necessary condition for optimality. If we further assume that the solution is optimal, then along the characteristic $\bm x (t; \bm x_0)$ we have that
\begin{equation}
\label{eq: control and value from BVP}
\bm u^* (t, \bm x) = \bm u^* (t; \bm x_0) ,
\qquad
V (t, \bm x) = v(t; \bm x_0) ,
\qquad
V_{\bm x} (t, \bm x) = \bm \lambda (t; \bm x_0) .
\end{equation}

In \cite{Kang2015, Kang2017_COA}, the two-point BVP \cref{eq: BVP} is solved for each point in a sparse grid. Applying \cref{eq: control and value from BVP}, the value function and its gradient are then calculated using high-dimensional interpolation. This technique is called the sparse grid characteristics method. But even in a sparse grid the number of points grows like $O \left( N (\log N)^{n-1} \right)$, where $n$ is the state dimension and $N$ is the number of grid points in each dimension. Thus one may have to solve a prohibitively large number of BVPs for higher-dimensional problems. Instead of sparse grid interpolation, we use data from solved BVPs to train a NN to approximate the value function. This approach is completely grid-free and hence applicable in high dimensions.

\begin{remark}
In general, the BVP admits multiple solutions which can sometimes be sub-optimal. The characteristics of the value function satisfy \cref{eq: BVP}, but there may be other solutions to these equations which are sub-optimal and therefore not characteristics of the value function. In many problems it is also possible for the characteristics to intersect, giving rise to non-smooth value functions and difficulties in applying \cref{eq: costate as gradient}.

Optimality of solutions to the BVP can be guaranteed under some convexity conditions (see e.g. \cite{Mangasarian1966}). For most dynamical systems it is difficult to verify such conditions globally, but we can guarantee optimality locally around an equilibrium point \cite{Lukes1969}. Addressing the challenge of global optimality in a broader context is beyond the scope of the present work, so in this paper we assume that solutions to the two-point BVP \cref{eq: BVP} are optimal. Under this assumption, the relationship between PMP and the value function as given in \cref{eq: control and value from BVP} holds everywhere.

Note the proposed method can still be applied to problems where this assumption cannot be verified. In such cases PMP remains the prevailing tool for computing candidate optimal solutions, and from these the proposed method will yield a feedback controller which satisfies necessary conditions for optimality.
\end{remark}

\subsection{Causality-free data generation}
\label{sec: time marching}

While solving the BVP is easier than solving the full HJB equation, we know of no general algorithm that is reliable and fast enough for real-time applications. However, in our approach the real-time feedback control computation is done by a NN which is trained \textit{offline}. Thus we can solve the BVP offline to generate data for training and evaluating such a NN. For this purpose, numerically solving the BVP can be manageable although it may require some parameter tuning. In this paper, we use an implementation of the BVP solver introduced in \cite{Kierzenka2001}. This algorithm is based on a three-stage Lobatto IIIa discretization, a collocation formula which provides a solution that is fourth-order accurate. But the algorithm is highly sensitive to the initial guess for $\bm x (t)$ and $\bm \lambda (t)$: there is no guarantee of convergence with an arbitrary initial guess, and in most cases a good initial guess for $\bm \lambda (t)$ cannot be derived from the problem physics. Furthermore, convergence is increasingly dependent on good initializations as we increase the length of the time interval.

To overcome this difficulty, we employ the \textit{time-marching} trick from \cite{Kang2015, Kang2017_COA}. This is a continuation technique in which we sequentially extend the solution from an initially short time interval to the final time $t_f$. Specifically, we choose a sequence of intermediate times
$$
0< t_1<t_2< \cdots < t_K=t_f,
$$
in which $t_1$ is small. For the short time interval $[0, t_1]$, the BVP solver converges given most initial guesses near the initial state $\bm x_0$. Then the resulting trajectory is rescaled over the longer time interval $[0, t_2]$. The rescaled trajectory is used as the initial guess to solve the BVP over $0\leq t \leq t_2$. We repeat this process until $t_K=t_f$, at which we obtain the full solution. By appropriately tuning the time sequence $\{t_k\}_{k=1}^K$, we can largely overcome the problem of sensitivity to initial guesses.

Computing many such solutions becomes expensive, which means that generating the large data sets necessary to train a NN can be difficult. With this in mind, we use the time-marching trick only to generate a small initial data set, and adaptively adding more points during training. The key to doing this efficiently is simulating the system dynamics using the partially-trained NN to close the loop. The closed-loop trajectory and predicted costate provide good guesses for the optimal state and costate, so that we can immediately solve \cref{eq: BVP} for all of $[0, t_f]$. Besides being more computationally efficient than time-marching, this approach also requires no parameter tuning. Details are presented in \cref{sec: adaptive sampling}, and numerical comparisons between this method and time-marching are given in \cref{sec: satellite adaptive sampling,sec: Burgers adaptive sampling}.

As an alternative to either of these two approaches, one could use backpropagation as suggested in \cite{Izzo2019}. However, this method does not allow one to choose initial conditions independently and so cannot be considered fully causality-free.

\section{Neural network approximation of the value function}
\label{sec: NN}

Neural networks have become a popular tool for modeling high-dimensional functions, since they are not dependent on discretizing the state space. In this paper, we apply NNs to approximate solutions of the HJB equation and evaluate the resulting feedback control in real-time. Specifically, we carry out the following steps:

\begin{enumerate}
\item \textit{Initial data generation:} We compute the value function, $V(t,\bm x)$, along trajectories $\bm x (t)$ from initial conditions chosen by Monte Carlo sampling. Data is generated by solving the BVP as discussed in \cref{sec: time marching}. In this initial data generation step, we require relatively few data points since more data can be added later at little computational cost.

\item \textit{Model training:} Given this data set, we train a NN to approximate the value function. Learning is guided by the underlying structure of the problem, specifically by asking the NN to satisfy Eq. \cref{eq: costate as gradient}. In doing so, we regularize the model and make efficient use out of small data sets.

\item \textit{Adaptive data generation:} In the initial training phase we only have a small data set, so the NN only roughly approximates the value function. We now expand the data set by generating data in regions where the value function is likely to be steep or complicated, and thus difficult to learn. Generating additional data is made efficient by good initial guesses obtained from NN-in-the-loop simulations of the system dynamics.

\item \textit{Model refinement and validation:} We continue training the model and increasing the size of the data set until we satisfy some convergence criteria. Then, we check the generalization accuracy of the trained NN on a new set of validation data computed at Monte Carlo sample points.

\item \textit{Feedback control:} We compute the feedback control online by evaluating the gradient of the trained NN and applying PMP. Notably, evaluation of the gradient is \textit{exact} and it is extremely cheap even for large $n$, enabling real-time implementation in high-dimensional systems.
\end{enumerate}

The crux of the proposed method depends on modeling the value function \cref{eq: value function} over a semi-global domain $\mathbb X \subset \mathbb R^n$. We present details of this process in the following subsections. In \cref{sec: intro NNs}, we review the basic structure of feedforward NNs and describe how we train a NN to model the value function. Then in \cref{sec: value gradient}, we propose a simple way to incorporate information about the known solution structure into training.  Finally in \cref{sec: NN in closed-loop system}, we demonstrate how to use the trained NN for feedback control. The adaptive data generation scheme is treated separately in \cref{sec: adaptive sampling}. The proposed method is illustrated in \cref{sec: rigid body} by solving a practical optimal attitude control problem for a rigid body satellite, and then applied to solve larger problems in \cref{sec: Burgers}.

\subsection{Feedforward neural networks}
\label{sec: intro NNs}

In this paper we use multilayer feedforward NNs. While many more sophisticated architectures have been developed for other applications, we find this basic architecture to be more than adequate for our purposes. Let $V (\cdot)$ be the function we wish to approximate and $V^{\text{NN}} (\cdot)$ be its NN representation. Feedforward NNs approximate complicated nonlinear functions by a composition of simpler functions, namely
$$
V (t, \bm x) \approx V^{\text{NN}} (t, \bm x)
	= g_L \circ \bm g_{L-1} \circ \cdots \circ \bm g_\ell \circ \cdots \circ \bm g_1 (t, \bm x) ,
$$
where each layer $\bm g_\ell (\cdot)$ is defined as
$$
\bm g_\ell (\bm y) = \sigma_\ell (\bm W_\ell \bm y + \bm b_\ell) .
$$
Here $\bm W_\ell$ and $\bm b_\ell$ are the weight matrices and bias vectors, respectively. $\sigma_\ell (\cdot)$ represents a nonlinear \textit{activation function} applied component-wise to its argument; popular choices include ReLU, tanh, and other similar functions. In this paper, we use tanh for all the hidden layers. The final layer, $g_L (\cdot)$, is typically linear, so $\sigma_L (\cdot)$ is the identity function.

Let $\bm \theta$ denote the collection of the parameters of the NN, i.e.
$$
\bm \theta \coloneqq \{ \bm W_\ell, \bm b_\ell \}_{\ell=1}^L .
$$
The NN is trained by optimizing over the parameters $\bm \theta$ to best approximate $V (t, \bm x)$ by $V^{\text{NN}} (t, \bm x; \bm \theta)$. Specifically, by solving the BVP \cref{eq: BVP} from a set of randomly sampled initial conditions, we get a data set
$$
\mathcal D = \left \{ \left( t^{(i)}, \bm x^{(i)} \right), V^{(i)} \right \}_{i=1}^{N_d} ,
$$
where $\left( t^{(i)}, \bm x^{(i)} \right)$ are the inputs, $V^{(i)} \coloneqq V \left( t^{(i)}, \bm x^{(i)} \right)$ are the outputs to be modeled, and $i = 1, 2, \dots, N_d$ are the indices of the data points. In the most na\"ive setting, the NN is then trained by solving the nonlinear regression problem,
\begin{equation}
\label{eq: regression problem}
\begin{array}{rl}
\underset{\bm \theta}{\text{minimize}} & \dfrac{1}{N_d} \displaystyle \sum_{i=1}^{N_d} \left[ V^{(i)} - V^{\text{NN}} \left(t^{(i)}, \bm x^{(i)} ; \bm \theta \right) \right]^2 .
\end{array}
\end{equation}

\subsection{Physics-informed machine learning}
\label{sec: value gradient}

Motivated by the development of physics-informed neural networks \cite{Raissi2019}, we expect that we can improve on the rudimentary loss function in \cref{eq: regression problem} by incorporating information about the underlying physics. In \cite{Raissi2019}, and in particular in the context of HJB equations in \cite{Khalaf2005,Cheng2007,Tassa2007,Sirignano2018}, the known underlying PDE and boundary conditions are imposed by minimizing a residual loss over spatio-temporal collocation points. In this approach, no data is gathered: the PDE is solved directly in the least-squares sense. But this residual must be evaluated over a large number of collocation points and can be rather expensive to compute. Thus we propose a simpler approach of modeling the costate $\bm \lambda (\cdot)$ along with the value function itself, taking full advantage of the ability to gather data along the characteristics of the HJB PDE.

Specifically, we know that the costate must satisfy Eq. \cref{eq: costate as gradient}, so we train the NN to minimize
$$
\left \Vert \bm \lambda (t; \bm x) - V^{\text{NN}}_{\bm x} (t, \bm x; \bm \theta) \right \Vert^2 ,
$$
where $V^{\text{NN}}_{\bm x} (\cdot)$ is the gradient of the NN model with respect to the state. This quantity is calculated using automatic differentiation. In machine learning, automatic differentiation is usually used to compute gradients with respect to the model parameters, but is just as easy to apply to computing gradients with respect to inputs. This gradient is exact, so no finite difference approximations are needed. In addition, the computational graph is pre-compiled so evaluating the gradient is cheap.

Costate data $\bm \lambda (t)$ is obtained for each trajectory as a natural product of solving the BVP \cref{eq: BVP}. Hence we have the augmented data set,
\begin{equation}
\label{eq: augmented data}
\mathcal D = \left \{ \left( t^{(i)}, \bm x^{(i)} \right), \left( V^{(i)}, \bm \lambda^{(i)} \right) \right \}_{i=1}^{N_d} ,
\end{equation}
where $\bm \lambda^{(i)} \coloneqq \bm \lambda \left( t^{(i)}; \bm x^{(i)} \right)$. We now define the {\em physics-informed learning problem},
\begin{equation}
\label{eq: value gradient problem}
\begin{array}{rl}
\underset{\bm \theta}{\text{minimize}}
	& \text{loss} \left( \bm \theta; \mathcal D \right)
	\coloneqq \underset{V}{\text{loss}} \left( \bm \theta; \mathcal D \right) + \mu \cdot \underset{\bm \lambda}{\text{loss}} \left( \bm \theta; \mathcal D \right) .
\end{array}
\end{equation}
Here $\mu \geq 0$ is a scalar weight, the loss with respect to data is
\begin{equation}
\label{eq: loss V}
\underset{V}{\text{loss}} \left( \bm \theta; \mathcal D \right)
	\coloneqq \frac{1}{N_d} \sum_{i=1}^{N_d} \left[ V^{(i)} - V^{\text{NN}} \left( t^{(i)}, \bm x^{(i)}; \bm \theta \right) \right]^2 ,
\end{equation}
and the gradient regularization is defined as
\begin{equation}
\label{eq: gradient loss}
\underset{\bm \lambda}{\text{loss}} \left( \bm \theta; \mathcal D \right)
	\coloneqq \frac{1}{N_d} \sum_{i=1}^{N_d} \left \Vert \bm \lambda^{(i)} - V^{\text{NN}}_{\bm x} \left( t^{(i)}, \bm x^{(i)} ; \bm \theta\right) \right \Vert^2 .
\end{equation}
Following standard practice, when computing the loss functions \cref{eq: loss V,eq: gradient loss}, the output data is linearly scaled to the range $[-1,1]$ to improve the scaling of the optimization problem.

A NN trained to minimize \cref{eq: value gradient problem} learns not just to fit the value data, but it is rewarded for doing so in a way that respects the underlying structure of the problem. Gradient regularization takes the known solution structure into account; this makes it preferable to the usual $L^1$ or $L^2$ regularization, which are based on the (heuristic) principle that simpler representations of data are likely to generalize better. Furthermore, we recall that the optimal control depends explicitly on $V_{\bm x} (\cdot)$ -- see Eqs. \cref{eq: optimal feedback control using dVdx} and \cref{eq: feedback control for control affine system}. Accurate approximation of $V_{\bm x} (\cdot)$ is therefore essential for calculating optimal controls. Our method achieves this through automatic differentiation to compute \textit{exact} gradients and by minimization of the gradient loss term \cref{eq: gradient loss}.

\subsection{Model validation}
\label{sec: validation}

In common practice, one randomly partitions the given data set \cref{eq: augmented data} into a training set $\mathcal D_{\text{train}}$ and validation set $\mathcal D_{\text{val}}$. During training, the loss functions \cref{eq: loss V} and \cref{eq: gradient loss} are calculated with respect to the training data $\mathcal D_{\text{train}}$. We then evaluate the performance of the NN against the validation data $\mathcal D_{\text{val}}$, which it did not observe during training. Good validation performance indicates that the NN generalizes well, i.e. it did not overfit the training data. We make the validation test more stringent by generating $\mathcal D_{\text{train}}$ and $\mathcal D_{\text{val}}$ from \textit{independently drawn} initial conditions, so that the two data sets do not share any part of the same trajectories.

We consider the following error metrics for validation. First, the relative mean absolute error (RMAE) of value function prediction, which is defined as
\begin{equation}
\text{RMAE} (\bm \theta; \mathcal D_{\text{val}}) \coloneqq \frac{\sum_{i=1}^{N_d} \left| V^{(i)} - V^{\text{NN}} \left( t^{(i)}, \bm x^{(i)} ; \bm \theta \right) \right|}{\sum_{i=1}^{N_d} \left| V^{(i)} \right|} .
\end{equation}
We also measure the relative mean $L^2$ error (RM$L^2$) of gradient prediction, which is defined as
\begin{equation}
\text{RM}L^2 (\bm \theta; \mathcal D_{\text{val}}) \coloneqq \frac{\sum_{i=1}^{N_d} \left \Vert \bm \lambda^{(i)} - V^{\text{NN}}_{\bm x} \left( t^{(i)}, \bm x^{(i)} ; \bm \theta \right) \right \Vert_2}{\sum_{i=1}^{N_d} \left \Vert \bm \lambda^{(i)} \right \Vert_2} .
\end{equation}
We consider these error metrics instead of pointwise relative errors in order to emphasize predictive accuracy in regions where a lot of control effort is needed. This is important because we are interested in designing nonlinear controllers which are effective and efficient far away from the equilibrium.

\subsection{Neural network in the closed-loop system}
\label{sec: NN in closed-loop system}

Once the NN is trained, evaluating $V^{\text{NN}}_{\bm x} (t, \bm x)$ at new inputs is highly efficient. Moreover, since we minimized the gradient loss \cref{eq: gradient loss} during training, we also expect $V^{\text{NN}}_{\bm x} (t, \bm x)$ to approximate the true gradient well. At runtime, whenever the feedback control needs to be computed, we evaluate $V^{\text{NN}}_{\bm x} \left( t, \bm x \right)$ and then solve \cref{eq: optimal feedback control using dVdx} based on this approximation.

For many problems of interest, the optimization problem \cref{eq: optimal feedback control using dVdx} admits an analytic or semi-analytic solution. In particular, for the important class of control affine systems with running cost convex in $\bm u$, we can solve \cref{eq: optimal feedback control using dVdx} analytically. Suppose that the system dynamics can be written in the form
$$
\dot{\bm x} = \bm f (t, \bm x) + \bm g (t, \bm x) \bm u ,
$$
where $\bm f (t, \bm x) : [0, t_f] \times \mathbb X \to \mathbb R^n$, $\bm g (t, \bm x) : [0, t_f] \times \mathbb X \to \mathbb R^{n \times m}$, and the control is unconstrained. Further, suppose that the running cost is of the form
$$
\mathcal L (t, \bm x, \bm u) = h (t, \bm x) + \bm u^T \bm {Wu} ,
$$
for some convex function $h (t, \bm x) : [0, t_f] \times \mathbb X \to \mathbb R$ and some positive definite weight matrix $\bm W \in \mathbb R^{m \times m}$. Then the Hamiltonian is
$$
\mathcal H (t, \bm x, \bm \lambda, \bm u) = h (t, \bm x) + \bm u^T \bm {Wu} + \bm \lambda^T \bm f (t, \bm x) + \bm \lambda^T \bm g (t, \bm x) \bm u .
$$
Now we apply PMP, which for unconstrained control requires
$$
\bm 0_{m \times 1} = \mathcal H_{\bm u} (t, \bm x, \bm \lambda, \bm u^*)
	= 2 \bm {Wu}^* + \bm g^T (t, \bm x) \bm \lambda .
$$
Letting $\bm \lambda = V_{\bm x} \left( t, \bm x \right)$ and solving for $\bm u^*$ yields the optimal feedback control law in explicit form:
\begin{equation}
\label{eq: feedback control for control affine system}
\bm u^* (t, \bm x; V_{\bm x}) = - \frac{1}{2} \bm W^{-1} \bm g^T (t, \bm x) V_{\bm x} (t, \bm x) .
\end{equation}
The resulting NN controller is then simply
\begin{equation}
\label{eq: NN feedback control}
\bm u^{\text{NN}} (t, \bm x) = \bm u^* \left( t, \bm x; V^{\text{NN}}_{\bm x} \right) = - \frac{1}{2} \bm W^{-1} \bm g^T (t, \bm x) V_{\bm x}^{\text{NN}} (t, \bm x) .
\end{equation}

\section{Adaptive sampling and model refinement}
\label{sec: adaptive sampling}

Since generating just a single data point requires solving a challenging BVP, it can be expensive to generate large data sets which adequately represent the value function. This necessitates training using limited data and a method to generate new data in a smart and efficient way. In this paper, effective training with small data sets is accomplished by incorporating information about the costate as discussed in \cref{sec: value gradient}, but also by combining progressive data generation with an efficient adaptive sampling technique.

Optimization methods in machine learning (see e.g. \cite{Bottou2018} for a comprehensive survey) are typically divided into second and first order methods. Second order methods like L-BFGS \cite{Byrd1995} rely on accurate gradient computations, and hence generally have to use the entire data set. For this reason they are often referred to as batch or full-batch methods. On the other hand, first order methods based on stochastic gradient descent (SGD) use only small subsets, or mini-batches, of the full data set. That is, at each optimization iteration $k$, the loss functions in \cref{eq: regression problem} and \cref{eq: value gradient problem} are evaluated only on a subset $\mathcal S_k \subset \mathcal D_{\text{train}}$ with $\left| \mathcal S_k \right| \ll \left| \mathcal D_{\text{train}} \right|$. Here $| \mathcal D |$ denotes the number of data points in a data set $\mathcal D$. Although second order methods converge much more quickly than first order methods, the necessary gradient calculations are prohibitively expensive for large data sets. Consequently, SGD variants have become the de facto standard for machine learning applications.

But in the context of deep learning, our NNs are small and data sets smaller. Thus we expect second order methods to be superior for our purposes. With a small initial data set, which we denote by $\mathcal D_{\text{train}}^1$, we find that training a low-fidelity model is very fast using L-BFGS. After this initial round, we want to increase the size of the data set so that it better captures the features of the value function. We then continue training the model using this larger data set, $\mathcal D_{\text{train}}^2$. We continue this process until some convergence conditions are satisfied.

Our approach is similar to and inspired by a progressive batching method proposed in \cite{Byrd2012}. The primary difference is that the problem addressed in \cite{Byrd2012} is a standard machine learning problem, where a massive data set is available from the start. This allows one to increase the sample size every few iterations, and take a completely different sample from the available data. In our problem, start with only a small amount of data and we can generate more as we go, but since data generation is expensive, we would like to generate only as much as is needed.

\subsection{Convergence test and sample size selection}

In this section we derive a convergence test and sample size selection scheme for the purpose of progressive data generation. To start, suppose that the internal optimizer (e.g. L-BFGS) converges in optimization round $r$ and let $\mathcal D_{\text{train}}^r$ be the available training data set. Given convergence of the internal optimizer, the first order necessary condition for optimality holds, so
\begin{equation}
\label{eq: standard first order optimality condition}
\left \Vert \frac{\del \text{loss}}{\del \bm \theta} \left( \bm \theta; \mathcal D_{\text{train}}^r \right) \right \Vert \ll 1 .
\end{equation}
Here $\text{loss} (\cdot)$ is the physics-informed loss defined in Eq. \cref{eq: value gradient problem}, and $\frac{\del \text{loss}}{\del \bm \theta} (\cdot)$ is its gradient with respect to the NN parameters $\bm \theta$. For true first order optimality, we would like the gradient to be small when evaluated over the entire continuous domain of interest, $[0, t_f] \times \mathbb X$. In other words, we want
\begin{equation}
\label{eq: ideal first order optimality condition}
\left \Vert \frac{\del \text{loss}}{\del \bm \theta} \left( \bm \theta; [0, t_f] \times \mathbb X \right) \right \Vert
\ll 1 ,
\end{equation}
where the Monte Carlo sums in Eqs. \cref{eq: loss V,eq: gradient loss} become integrals in the limit as the size of the data set approaches infinity.

The simplest way to see if \cref{eq: ideal first order optimality condition} holds is to generate a validation data set $\mathcal D_{\text{val}}$. Then using the fact that $\frac{\del \text{loss}}{\del \bm \theta} \left( \bm \theta; \mathcal D_{\text{val}} \right) \to \frac{\del \text{loss}}{\del \bm \theta} \left( \bm \theta; [0, t_f] \times \mathbb X \right)$ in the limit as $\left| \mathcal D_{\text{val}} \right| \to \infty$, one checks if, for example,
\begin{equation}
\label{eq: validation test}
\left \Vert \frac{\del \text{loss}}{\del \bm \theta} \left( \bm \theta; \mathcal D_{\text{val}} \right) \right \Vert < \epsilon ,
\end{equation}
for some small parameter $\epsilon > 0$. Convergence tests like \cref{eq: validation test} are standard in machine learning and are useful for testing generalization performance. But for many practical problems, it may be too expensive to generate enough validation data to make the test meaningful. More importantly, such tests provides no clear guidance in selecting the sample size $\left| \mathcal D_{\text{train}}^{r+1} \right|$ should they not be satisfied.

In this paper, we use validation tests to quantify model accuracy after training is complete (see \cref{sec: validation}). Indeed, the ability to empirically validate solutions is a key benefit of the causality-free approach. For the purpose of determining convergence between training rounds, however, we propose a different statistically motivated test which provides information on choosing $\left| \mathcal D_{\text{train}}^{r+1} \right|$. The idea is simple: since we already assume \cref{eq: standard first order optimality condition} holds, then to ensure that \cref{eq: ideal first order optimality condition} is also satisfied, it suffices to check that the error in approximating \cref{eq: ideal first order optimality condition} by \cref{eq: standard first order optimality condition} is relatively small.

To motivate this more rigorously, consider a finite sample set $\mathcal D \subset [0, t_f] \times \mathbb X$ with fixed size $| \mathcal D |$, and assume that the sample points $\left( t^{(i)}, \bm x^{(i)} \right) \in \mathcal D$ are independent and identically distributed (i.i.d.)\footnote{In practice, while initial conditions are i.i.d., points at future times lie along the optimal trajectories coming from these initial conditions and are thus spatially correlated. Adaptive sampling (see \cref{sec: adaptive sampling implementation}) also introduces sample dependence. This likely reduces sample variance compared to i.i.d. data, but we still find the numerical tests useful for providing sample size guidelines. In addition, if we learn only the initial-time value function as in \cref{sec: rigid body}, then sample independence can be upheld if we forego adaptive sample placement.}. By design, if $\left( t^{(i)}, \bm x^{(i)} \right)$ are i.i.d. then the sample gradient $\frac{\del \text{loss}}{\del \bm \theta} \left( \bm \theta ; \mathcal D \right)$ is an unbiased estimator for the true gradient (evaluated over the entire continuous domain). That is,
\begin{equation}
\label{eq: expectation of gradient}
\mathbb E_{\mathcal D} \left[ \frac{\del \text{loss}}{\del \bm \theta} \left( \bm \theta; \mathcal D \right) \right]
	= \frac{\del \text{loss}}{\del \bm \theta} \left( \bm \theta; [0, t_f] \times \mathbb X \right) ,
\end{equation}
where $\mathbb E_{\mathcal D} [\cdot] \coloneqq \mathbb E_{\mathcal D \subset [0, t_f] \times \mathbb X} [ \cdot ]$ denotes the population mean over all possible finite sample sets $\mathcal D \subset [0, t_f] \times \mathbb X$ with fixed size $| \mathcal D |$. Intuitively, \cref{eq: expectation of gradient} implies that if \cref{eq: standard first order optimality condition} holds, then on average we also have \cref{eq: ideal first order optimality condition}, as desired. But we must control the mean square error (MSE) of the estimator, which is given by
\begin{align*}
\numberthis
\text{MSE} \left[ \frac{\del \text{loss}}{\del \bm \theta} \left( \bm \theta; \mathcal D \right) \right]
	\coloneqq& \mathbb E_{\mathcal D} \left[ \left \Vert \frac{\del \text{loss}}{\del \bm \theta} \left( \bm \theta; \mathcal D \right) - \frac{\del \text{loss}}{\del \bm \theta} \left( \bm \theta; [t_0, t_f] \times \mathbb X \right) \right \Vert^2 \right] \\
	=& \mathbb E_{\mathcal D} \left[ \sum_{j=1}^{| \bm \theta |} \left( \frac{\del \text{loss}}{\del \theta_j} (\bm \theta; \mathcal D) - \frac{\del \text{loss}}{\del \theta_j} (\bm \theta ; [t_0, t_f] \times \mathbb X) \right)^2 \right] .
\end{align*}
To simplify this, using linearity of the expectation we obtain
$$
\text{MSE} \left[ \frac{\del \text{loss}}{\del \bm \theta} \left( \bm \theta; \mathcal D \right) \right]
	= \sum_{j=1}^{| \bm \theta |} \text{Var}_{\mathcal D \subset [0, t_f] \times \mathbb X} \left[ \frac{\del \text{loss}}{\del \theta_j} (\bm \theta; \mathcal D) \right] ,
$$
and then by construction of the loss function,
$$
\text{MSE} \left[ \frac{\del \text{loss}}{\del \bm \theta} \left( \bm \theta; \mathcal D \right) \right]
	= \sum_{j=1}^{| \bm \theta |} \text{Var}_{\mathcal D \subset [0, t_f] \times \mathbb X} \left[ \frac{1}{\left| \mathcal D \right|} \sum_{i=1}^{\left| \mathcal D \right|} \frac{\del \text{loss}}{\del \theta_j} \left( \bm \theta; \left( t^{(i)}, \bm x^{(i)} \right) \right) \right] .
$$
Using the simplifying assumption that $\left( t^{(i)}, \bm x^{(i)} \right)$ are i.i.d., this becomes
\begin{align*}
\text{MSE} \left[ \frac{\del \text{loss}}{\del \bm \theta} \left( \bm \theta; \mathcal D \right) \right]
	=& \frac{1}{\left| \mathcal D \right|^2} \sum_{j=1}^{| \bm \theta |} \sum_{i=1}^{\left| \mathcal D \right|} \text{Var}_{(t, \bm x) \in [0, t_f] \times \mathbb X} \left[ \frac{\del \text{loss}}{\del \theta_j} \left( \bm \theta; ( t, \bm x) \right) \right] \\
\numberthis
\label{eq: estimator variance}
	=& \frac{1}{\left| \mathcal D \right|} \sum_{j=1}^{| \bm \theta |} \text{Var}_{(t, \bm x) \in [0, t_f] \times \mathbb X} \left[ \frac{\del \text{loss}}{\del \theta_j} \left( \bm \theta; (t, \bm x) \right) \right] .
\end{align*}
If the estimation error is small, then the sample mean is likely to be a good approximation of the true mean. Hence we expect that $\left \Vert \frac{\del \text{loss}}{\del \bm \theta} \left( \bm \theta; [0, t_f] \times \mathbb X \right) \right \Vert$ will also be small as desired. To this end, we require that the root MSE not be too large compared to the expected gradient. Specifically, we check if
\begin{equation}
\label{eq: norm test base}
\sqrt{ \text{MSE} \left[ \frac{\del \text{loss}}{\del \bm \theta} \left( \bm \theta; \mathcal D \right) \right] }
	\leq C \left \Vert \mathbb E_{\mathcal D} \left[ \frac{\del \text{loss}}{\del \bm \theta} \left( \bm \theta; \mathcal D \right) \right] \right \Vert_1 ,
\end{equation}
where $C > 0$ is a scalar parameter. On the right hand side we use the $L^1$ norm instead of the $L^2$ as it is less sensitive to outliers in the loss gradient. In practice we find that this makes the test less likely to suggest unreasonably large sample sizes.

In practice, evaluating of the true population variances on the left hand side of \cref{eq: norm test base} is computationally intractable. But we can approximate these by the corresponding \textit{sample} variances\footnote{Computing a large number of individual gradients can still be too costly, so we often evaluate sample variances over a smaller subset of the training data.} taken over all data $\left( t^{(i)}, \bm x^{(i)} \right) \in \mathcal D_{\text{train}}^r$, which we denote by $\text{Var}_{\mathcal D_{\text{train}}^r} [\cdot] \coloneqq \text{Var}_{\left( t^{(i)}, \bm x^{(i)} \right) \in \mathcal D_{\text{train}}^r} [\cdot] $:
$$
\text{MSE} \left[ \frac{\del \text{loss}}{\del \bm \theta} \left( \bm \theta; \mathcal D \right) \right]
	\approx \frac{1}{\left| \mathcal D_{\text{train}}^r \right|} \sum_{j=1}^{| \bm \theta |} \text{Var}_{\mathcal D_{\text{train}}^r} \left[ \frac{\del \text{loss}}{\del \theta_j} \left( \bm \theta; \left( t^{(i)}, \bm x^{(i)} \right) \right) \right] .
$$
Similarly, we approximate the expected gradient on the right hand side of \cref{eq: norm test base} by the sample gradient and arrive at the following practical convergence criterion:
\begin{equation}
\label{eq: our norm test}
\sqrt{\sum_{j=1}^{| \bm \theta |} \text{Var}_{\mathcal D_{\text{train}}^r} \left[ \frac{\del \text{loss}}{\del \theta_j} \left( \bm \theta; \left( t^{(i)}, \bm x^{(i)} \right) \right) \right] }
	\leq C \left \Vert \frac{\del \text{loss}}{\del \bm \theta} \left( \bm \theta; \mathcal D_{\text{train}}^r \right) \right \Vert_1 \sqrt{ \left| \mathcal D_{\text{train}}^r \right| } .
\end{equation}

If the convergence test \cref{eq: our norm test} is satisfied, then it is likely that the expected gradient $\left \Vert \frac{\del \text{loss}}{\del \bm \theta} \left( \bm \theta; [0, t_f] \times \mathbb X \right) \right \Vert$ is also small. In other words, we expect that the parameters $\bm \theta$ satisfies the first order optimality conditions evaluated over the entire domain, so we can stop optimization. Satisfaction of \cref{eq: our norm test} does not imply that the trained model is good -- merely that seeing more data would probably not improve it significantly. On the other hand, when the criterion is not met, then it guides us in selecting the next sample size $\left| \mathcal D_{\text{train}}^{r+1} \right|$. Concretely, suppose that the ratio of the sample variance to the sample gradient doesn't change significantly by increasing the size of the data set, i.e.
$$
\frac{\sqrt{\sum_{j=1}^{| \bm \theta |} \text{Var}_{\mathcal D_{\text{train}}^{r+1}} \left[ \frac{\del  \text{loss}}{\del \theta_j} \left( \bm \theta; \left( t^{(i)}, \bm x^{(i)} \right) \right) \right]}}
{\left \Vert \frac{\del \text{loss}}{\del \bm \theta} \left( \bm \theta; \mathcal D_{\text{train}}^{r+1} \right) \right \Vert_1}
\approx
\frac{\sqrt{\sum_{j=1}^{| \bm \theta |} \text{Var}_{\mathcal D_{\text{train}}^r} \left[ \frac{\del  \text{loss}}{\del \theta_j} \left( \bm \theta; \left( t^{(i)}, \bm x^{(i)} \right) \right) \right]}}
{\left \Vert \frac{\del \text{loss}}{\del \bm \theta} \left( \bm \theta; \mathcal D_{\text{train}}^r \right) \right \Vert_1} .
$$
Then the appropriate choice of $\left| \mathcal D_{\text{train}}^{r+1} \right|$ to satisfy \cref{eq: our norm test} after the next round is such that
\begin{equation}
\label{eq: our sample size selection}
M \left| \mathcal D_{\text{train}}^r \right| \geq \left| \mathcal D_{\text{train}}^{r+1} \right| \geq
	\frac{\sum_{j=1}^{| \bm \theta |} \text{Var}_{\mathcal D_{\text{train}}^r} \left[ \frac{\del \text{loss}}{\del \theta_j} \left( \bm \theta; \left( t^{(i)}, \bm x^{(i)} \right) \right) \right]}
{\left( C \left \Vert \frac{\del \text{loss}}{\del \bm \theta} \left( \bm \theta; \mathcal D_{\text{train}}^r \right) \right \Vert_1 \right)^2} ,
\end{equation}
where $M > 1$ is a scalar parameter which prevents the data set size from growing too quickly. Throughout this paper we use $M = 2$.

The convergence test \cref{eq: our norm test} and sample size selection scheme \cref{eq: our sample size selection} derived above are close to that used in \cite{Byrd2012}, except that we employ the $L^1$ norm of the sample gradient in the denominator instead of the $L^2$ norm. We prefer the $L^1$ norm because it is less sensitive to outliers in the loss gradient. Intuitively, this improves robustness by making the test less likely to suggest unreasonably large sample sizes. We also contribute a different derivation, coming from the perspective of progressive data generation as opposed to sampling from a large pre-existing data set. Finally, like \cite{Byrd2012} our results are not specific to learning solutions to the HJB equation. They can be applied to many data-driven optimization problems where data is scarce but can be generated over time. Notably, these results facilitate the use of existing algorithms for second order and constrained optimization in such applications.

\subsection{Adaptive data generation with NN warm start}
\label{sec: adaptive sampling implementation}

The sample size selection criterion \cref{eq: our sample size selection} we propose indicates how many data are necessary to satisfy the convergence test \cref{eq: our norm test}, assuming a uniform sampling from the domain. In practice, since all the data we generate will be new, we can choose to generate new data where it is needed most, hence the term \textit{adaptive sampling}. This condition for generating new data can be interpreted in many ways. In this paper, we concentrate samples where $\left \Vert V^{\text{NN}}_{\bm x} (\cdot) \right \Vert$ is large. Regions of the value function with large gradients tend to be steep or complicated, and thus may benefit from having more data to learn from. Furthermore, these regions correspond to places where the control effort is large and hence we would like controllers to be especially accurate there.

Specifically, for each initial condition we want to integrate, we can first randomly sample a set of $N_c$ candidate initial conditions from $\mathbb X$. A quick pass through the NN yields the predicted gradient at all candidate points:
$$
\left \{ V^{\text{NN}}_{\bm x}  \left( 0, \bm x_0^{(i)} \right) \right \}_{i=1}^{N_c} .
$$
We then choose the point(s) with the largest predicted gradient norms and solve the BVP \cref{eq: BVP} for each of these. To aid in solving these BVPs, instead of using the time-marching trick described in \cref{sec: time marching}, we simulate the system dynamics using the partially-trained NN as the closed-loop controller and predicting $\bm \lambda (t) \approx V^{\text{NN}}_{\bm x} (t, \bm x (t))$ along the trajectory. In most cases, this yields an approximate solution which is reasonably close to the optimal state and costate. By supplying this trajectory as an initial guess to the BVP solver, we then quickly and reliably obtain a solution to the BVP for the full time interval $[0, t_f]$. This process is repeated for new initial conditions until we obtain the desired amount of data (each trajectory may contain hundreds of data points). We refer to this technique as a \textit{NN warm start}. A summary of the full training procedure is given in \cref{alg: training}.

\Cref{alg: training} enables us to build up a rich data set and a high-fidelity model of $V (\cdot)$. Moreover, the data set is not constrained to lie within a small neighborhood of some nominal trajectory. It can contain points from the entire domain $\mathbb X$, and we can concentrate more data near complicated features of the value function. As we progressively refine the NN model, we can adjust the gradient loss weight $\mu$, as well as other hyperparameters such as the internal optimizer convergence tolerance and the number of terms in the L-BFGS Hessian approximation. As the NN is already partially-trained, fewer iterations should be needed for convergence in each round so we can afford to make each iteration more expensive.

\begin{algorithm}[t!]
\caption{Adaptive sampling and model refinement}
\label{alg: training}
\begin{algorithmic}[1]
\STATE{Generate $\mathcal D_{\text{train}}^1$ using time-marching}
\FOR{$r = 1, 2, \dots$}
	\STATE{Solve \cref{eq: value gradient problem} for $\bm \theta$}
	\IF{\cref{eq: our norm test} is satisfied}
		\RETURN{optimized parameters $\bm \theta$ and NN validation accuracy}
	\ELSE
		\WHILE{\cref{eq: our sample size selection} is not satisfied}
			\STATE{Sample candidate initial conditions $\bm x_0^{(i)}$, $i = 1, \dots, N_c$}
			\STATE{In parallel, predict $\left \Vert V^{\text{NN}}_{\bm x}  \left( 0, \bm x_0^{(i)} \right) \right \Vert$, $i=1, \dots, N_c$}
			\STATE{Choose the initial condition(s) with largest predicted gradient norm and use NN warm start to solve the corresponding BVP(s) \cref{eq: BVP}}
			\STATE{Add the resulting trajectorie(s) to $\mathcal D_{\text{train}}^{r+1}$}
		\ENDWHILE
	\ENDIF
\ENDFOR
\end{algorithmic}
\end{algorithm}

\section{Application to rigid body attitude control}
\label{sec: rigid body}

To illustrate the capabilities of proposed method, we consider the six-state rigid body model of a satellite studied by Kang and Wilcox \cite{Kang2015, Kang2017_COA}. With the sparse grid characteristics method, they interpolate the value function at initial time, $V(t=0, \bm x)$, and use this for moving horizon feedback control of the nonlinear system. We use their successful results as a baseline for evaluating our method.

Let $\{ \bm e_1, \bm e_2, \bm e_3 \}$ be an inertial frame of orthonormal vectors and let $\{ \bm e_1', \bm e_2', \bm e_3' \}$ be a body frame. The state of the satellite is then written as $\bm x = \begin{pmatrix}
\bm v & \bm \omega
\end{pmatrix}$. Here $\bm v$ is the attitude of the satellite represented in Euler angles,
$$
\bm v = \begin{pmatrix}
\phi & \theta & \psi
\end{pmatrix}^T ,
$$
in which $\phi$, $\theta$, and $\psi$ are the angles of rotation around $\bm e_1'$, $\bm e_2'$, and $\bm e_3'$, respectively, in the order $(1, 2, 3)$. These are also commonly called roll, pitch, and yaw. $\bm \omega$ denotes the angular velocity in the body frame,
$$
\bm \omega = \begin{pmatrix}
\omega_1 & \omega_2 & \omega_3
\end{pmatrix}^T .
$$
For details see \cite{Diebel2006}. The state dynamics are
$$
\begin{pmatrix}
\dot{\bm v} \\
\bm J \dot{\bm \omega}
\end{pmatrix}
= \begin{pmatrix}
\bm E (\bm v) \bm \omega \\
\bm S (\bm \omega) \bm R (\bm v) \bm h + \bm B \bm u
\end{pmatrix} .
$$
Here $\bm E (\bm v) , \bm S (\bm \omega), \bm R (\bm v) : \mathbb R^3 \to \mathbb R^{3 \times 3}$ are matrix-valued functions defined as
$$
\bm E (\bm v) \coloneqq
\begin{pmatrix}
1 & \sin \phi \tan \theta & \cos \phi \tan \theta \\
0 & \cos \phi & - \sin \phi \\
0 & \sin \phi / \cos \theta & \cos \phi / \cos \theta
\end{pmatrix} ,
\qquad
\bm S (\bm \omega) \coloneqq
\begin{pmatrix}
0 & \omega_3 & - \omega_2 \\
- \omega_3 & 0 & \omega_1 \\
\omega_2 & - \omega_1 & 0
\end{pmatrix} ,
$$
and
$$
\bm R (\bm v) \coloneqq
\begin{pmatrix}
\cos \theta \cos \psi & \cos \theta \sin \psi & - \sin \theta \\
\sin \phi \sin \theta \cos \psi - \cos \phi \sin \psi & \sin \phi \sin \theta \sin \psi + \cos \phi \cos \psi & \cos \theta \sin \phi \\
\cos \phi \sin \theta \cos \psi + \sin \phi \sin \psi & \cos \phi \sin \theta \sin \psi - \sin \phi \cos \psi & \cos \theta \cos \phi
\end{pmatrix} .
$$
Further, $\bm J \in \mathbb R^{3 \times 3}$ is a combination of the inertia matrices of the momentum wheels and the rigid body without wheels, $\bm h \in \mathbb R^3$ is the total constant angular momentum of the system, and $\bm B \in \mathbb R^{3 \times m}$ is a constant matrix where $m$ is the number of momentum wheels. To control the system, we apply a torque $\bm u (t, \bm v, \bm \omega) : [0, t_f] \times \mathbb R^3 \times \mathbb R^3 \to \mathbb R^m$.

We consider the fully-actuated case where $m = 3$. Let
$$
\bm B = \begin{pmatrix}
1 & 1/20 & 1/10 \\
1/15 & 1 & 1/10 \\
1/10 & 1/15 & 1
\end{pmatrix} ,
\qquad
\bm J = \begin{pmatrix}
2 & 0 & 0 \\
0 & 3 & 0 \\
0 & 0 & 4
\end{pmatrix} ,
\qquad
\bm h = \begin{pmatrix}
1 \\
1 \\
1
\end{pmatrix} .
$$
The optimal control problem is
\begin{equation}
\label{eq: satellite OCP}
\left \{
\begin{array}{cl}
\underset{\bm u (\cdot)}{\text{minimize}} & J \left[ \bm u (\cdot) \right] = \displaystyle \int_t^{t_f} \mathcal L (\bm v, \bm \omega, \bm u) d\tau + \dfrac{W_4}{2} \Vert \bm v (t_f) \Vert^2 + \dfrac{W_5}{2} \Vert \bm \omega (t_f) \Vert^2 , \\
\text{subject to} & \dot{\bm v} = \bm E (\bm v) \bm \omega , \\
	& \bm J \dot{\bm \omega} = \bm S (\bm \omega) \bm R (\bm v) \bm h + \bm B \bm u .
\end{array}
\right .
\end{equation}
Here
$$
\mathcal L (\bm v, \bm \omega, \bm u) = \frac{W_1}{2} \Vert \bm v \Vert^2 + \frac{W_2}{2} \Vert \bm \omega \Vert^2 + \frac{W_3}{2} \Vert \bm u \Vert^2
$$
and
$$
\begin{array}{cccccc}
W_1 = 1,&
W_2 = 10,&
W_3 = \dfrac{1}{2},&
W_4 = 1,&
W_5 = 1,&
t_f = 20.
\end{array}
$$
Finally, we consider initial conditions in the domain
\begin{equation}
\label{eq: satellite initial conditions}
\mathbb X_0 = \left \{ \left. \bm v, \bm \omega \in \mathbb R^3 \right| - \frac{\pi}{3} \leq \phi, \theta, \psi \leq \frac{\pi}{3}
\text{ and }
- \frac{\pi}{4} \leq \omega_1, \omega_2, \omega_3 \leq \frac{\pi}{4} \right \} .
\end{equation}

In \cite{Kang2017_COA}, to avoid discretizing time the value function is approximated only at initial time $t = 0$. In order to facilitate a fair comparison we do the same. This means that we model $V (0, \bm v, \bm \omega) \approx V^{\text{NN}} (\bm v, \bm \omega)$, i.e. the NN does {\em not} take time as an input variable. Consequently the control is implemented with a time-independent moving horizon rather than as a time-dependent optimal control. In other words, at each time $t$ when we evaluate the control, we assume $t = 0$ and return $\bm u (t) = \bm u^{\text{NN}} \left( \bm v (t), \bm \omega (t) \right)$. Controlling the system using moving horizon feedback is standard practice. It is also reasonable for the present case because the problem dynamics are time-invariant and the time horizon is relatively long. Because of this we observe near-optimal performance from the moving horizon controller.

\subsection{Learning the value function}
\label{sec: satellite value function learning}

\begin{figure}[t!]
\centering
\includegraphics[width = .85 \textwidth]{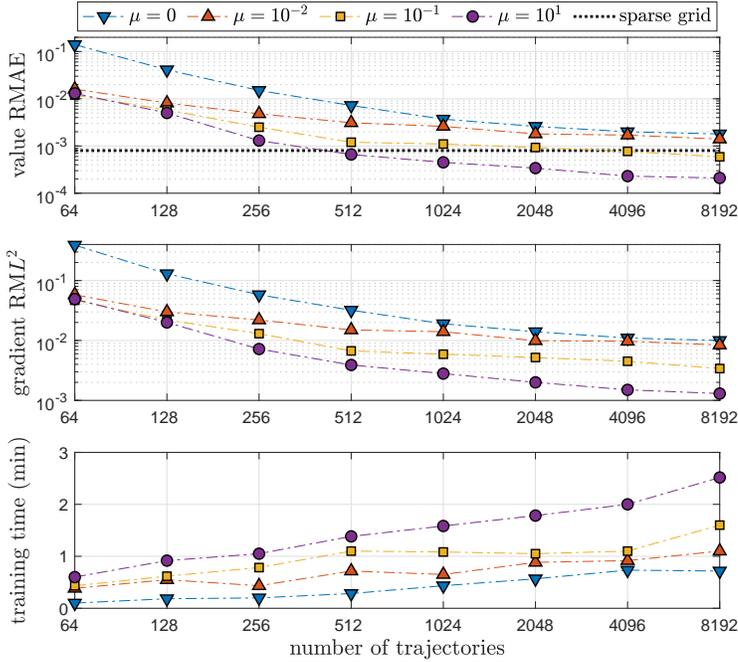}
\caption{Validation accuracy and training time of NNs for modeling the initial time value function $V (0, \bm v, \bm \omega)$ of the optimal attitude control problem \cref{eq: satellite OCP}. All NNs have the same parameter initialization and are run on an NVIDIA RTX 2080Ti GPU.}
\label{fig: satellite results t0}
\end{figure}

In this section, we present numerical results of our implementation of a NN for modeling the initial-time value function of the rigid body attitude control problem \cref{eq: satellite OCP}. To obtain data, we uniformly sample initial conditions $\left( \bm v^{(i)}, \bm \omega^{(i)} \right)$ from the domain $\mathbb X_0$ defined in \cref{eq: satellite initial conditions}, and for each initial value, solve the two-point BVP \cref{eq: BVP} using time-marching and the SciPy \cite{SciPy} implementation of the three-stage Lobatto IIIa algorithm in \cite{Kierzenka2001}. Each integrated trajectory contains around 100 data points on average, but we use only initial time data, $V \left( 0, \bm v^{(i)}, \bm \omega^{(i)} \right)$. For validation, we generate a data set containing $\left| \mathcal D_{\text{val}} \right| = 1000$ data points (at $t=0$), and keep this fixed throughout all the tests. As a baseline, the sparse grid characteristics method with $\left| G_{\text{sparse}}^{13} \right| = 44,698$ grid points achieves a RMAE of $8.00 \times 10^{-4}$ on this validation data set.

We implement a standard feedforward NN in TensorFlow 1.11 \cite{TensorFlow} and train it to approximate $V(0, \bm v, \bm \omega)$. The NN has three hidden layers with 64 neurons in each, but many alternate configurations of depth and width also work. For optimization, we use the SciPy interface for the L-BFGS optimizer \cite{SciPy, Byrd1995}. \Cref{fig: satellite results t0} displays the results of a series of tests in which we vary the weight $\mu$ on the value gradient loss term \cref{eq: gradient loss} and the size of the training data set. Results are compared to those obtained in \cite{Kang2017_COA}.

We highlight that with just 512 data points, we can train NNs with better accuracy than the sparse grid characteristics method with $\left| G_{\text{sparse}}^{13} \right| = 44,698$ points. Thus for this problem, the proposed method is about 90 times as data-efficient. With 8192 data points, the NN can be almost four times as accurate as the sparse grid characteristics method. This level of accuracy with small data sets is obtained only with physics-informed learning. In particular, NNs trained by pure regression \cref{eq: regression problem} cannot match the accuracy of the sparse grid characteristics method, as shown in \cref{fig: satellite results t0} for the case with $\mu=0$. Accuracy improves as we increase $\mu$ but with diminishing returns for $\mu \geq 10$. While physics-informed learning is more costly, it facilitates the use of much smaller data sets, and the increased training time is still quite short.

\subsection{Training with adaptive data generation}
\label{sec: satellite adaptive sampling}

\begin{figure}[t!]
\centering
\includegraphics[width = 0.85 \textwidth]{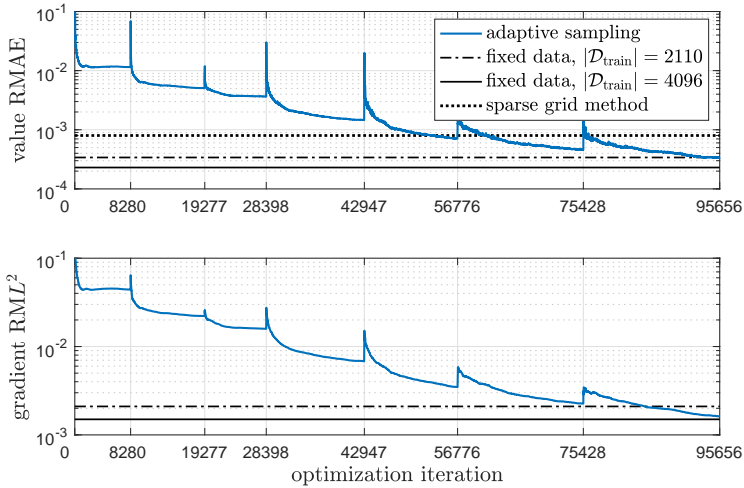}
\caption{Progress of adaptive sampling and model refinement for the rigid body problem \cref{eq: satellite OCP}, compared to training on fixed data sets and the sparse grid characteristics method. Spikes in the error correspond to the start of new training rounds and expansion of the training data set.}
\label{fig: adaptive training}
\end{figure}

Performing a thorough systematic study of the adaptive sampling and model refinement technique proposed in \cref{sec: adaptive sampling} is rather complicated, since a successful implementation depends on various hyperparameter settings, which can and perhaps should change each optimization round. Results also depend on random chance, since data points are chosen in a (partially) random way and the randomly-initialized NN training problem is highly non-convex. For this reason, in this section we show a just few conservative results which we feel illustrate the potential of the method.

\Cref{fig: adaptive training} shows the progress of the validation error during training when using adaptive sampling starting from a data set with $\left| \mathcal D_{\text{train}}^1 \right| = 64$ points. We set the gradient loss weight to $\mu = 10$ and the convergence parameter in \cref{eq: our norm test} to $C = 0.25$. After each round, we check the convergence criterion \cref{eq: our norm test} and increase the number of training data according to \cref{eq: our sample size selection}. Each data set includes all previously generated data, and we generate additional data as needed through \cref{alg: training}. With these configurations, the model passes the convergence test after seven training rounds and observing a total of $\left| \mathcal D_{\text{train}}^{7} \right| = 2110$ samples.

The final value function accuracy is $3.3 \times 10^{-4}$: over twice as accurate as the sparse grid method with about twenty times fewer data, and the gradient prediction accuracy is $1.6 \times 10^{-3}$. As shown in \cref{fig: adaptive training}, the gradient predictions of the network trained using the adaptive algorithm are just as accurate as a network trained on a fixed data set of $\left| \mathcal D_{\text{train}} \right| = 4096$ samples. That is to say, the adaptive sampling method facilitates more acurate gradient predictions using fewer data. These results highlight the main advantages of the adaptive sampling and model refinement method: the ability to overcome an initial lack of data, efficiently generate a large data set, and improve gradient prediction accuracy which is needed for effective control. To fully realize the potential of the method, hyperparameters like $\mu$, $C$, and internal optimizer parameters need to be adjusted in each round. Development of algorithms to do this adaptively remains a topic for future research.

\begin{table}[t!]
\centering
\begin{tabular}{| c | c | c |}
\hline
&&\\[-1em]
$K$ & $\%$ BVP convergence & mean integration time \\
\hline
&&\\[-1em]
1 & $0.3 \%$ & 0.37 s \\
2 & $38.7 \%$ & 0.44 s \\
3 & $76.2 \%$ & 0.40 s \\
4 & $92.9 \%$ & 0.45 s \\
8 & $98.4 \%$ & 0.53 s \\
\hline
\end{tabular}
\caption{Convergence of BVP solutions for \cref{eq: satellite OCP} when using the time-marching trick, depending on the number of steps in the sequence $\{ t_k \}_{k=1}^K$. The case $K=1$ corresponds to a direct solution attempt over the whole time interval with no time-marching. BVP integration time is measured only on successful attempts -- failed solution attempts usually take much longer.}
\label{table: satellite BVP convergence}
\begin{tabular}{| c | c | c | c | c |}
\hline
&&&&\\[-1em]
$\mu$ & training time & gradient RM$L^2$ & $\%$ BVP conv. & mean int. time \\
\hline
&&&&\\[-1em]
$10^{-8}$ & 7 s & $2.5 \times 10^{-1}$ & $88.0 \%$ & 0.50 s \\
$10^{-4}$ & 19 s & $1.4 \times 10^{-1}$ & $98.6 \%$ & 0.48 s \\
1 & 23 s & $4.5 \times 10^{-2}$ & $99.7\%$ & 0.44 s \\
\hline
\end{tabular}
\caption{Convergence of BVP solutions for \cref{eq: satellite OCP} when using NN warm start with NNs of varying gradient prediction accuracy. BVP integration time is measured only on successful attempts.}
\label{table: satellite BVP convergence NN}
\end{table}

Next, we investigate the convergence of the BVP solver with time-marching and NN warm start. Results are given in \cref{table: satellite BVP convergence} and \cref{table: satellite BVP convergence NN}, respectively. For these tests, we randomly sample $10^6$ candidate initial conditions and pick 1000 points with the largest predicted gradient norm, $\left \Vert V_{\bm x}^{\text{NN}} (\cdot) \right \Vert$. Initial conditions with large gradient norm tend to be located in regions where the value function is steep and the control effort is large, and may thus be more difficult to solve. The set of initial conditions is fixed for all tests.

In the first row of \cref{table: satellite BVP convergence}, we attempt to solve the BVP with no time-marching, i.e. over the entire time interval without constructing any initial guess. In this case, the proportion of convergent solutions is extremely small, obviating the need for good initial guesses. As shown in \cref{table: satellite BVP convergence}, we reliably obtain solutions for this problem when we use at least $K=4$ time intervals. We note that the initial conditions are purposefully chosen to be difficult -- if we simply take uniform samples from the domain $\mathbb X_0$, the proportion of convergent solutions increases significantly.

In \cref{table: satellite BVP convergence NN}, we present results using NN warm start. We train several NNs on a data set of only 64 points. Because the data set is so small, {\em each NN takes only seconds to finish training}. We also experiment with using different gradient loss weights $\mu$ for each NN. This directly impacts the accuracy in predicting the initial-time costate, $\bm \lambda (0; \bm v_0, \bm \omega_0) \approx V^{\text{NN}}_{\bm x} ( \bm v_0, \bm \omega_0)$, which in turn is key to synthezing optimal controls.

Even with these low-fidelity models, the rate of BVP convergence is just as high as when using $K=4$ time intervals for time-marching. The quality of initial guesses improves with better costate prediction, and it is not difficult to exceed $99 \%$ convergence. For this problem, the speed of the two methods is about the same. However, when we consider higher-dimensional problems in \cref{sec: Burgers adaptive sampling}, we find that NN warm start significantly improves both reliability and efficiency.

\subsection{Closed-loop simulation}

\begin{figure}[t!]
\centering
\includegraphics[width = .85 \textwidth]{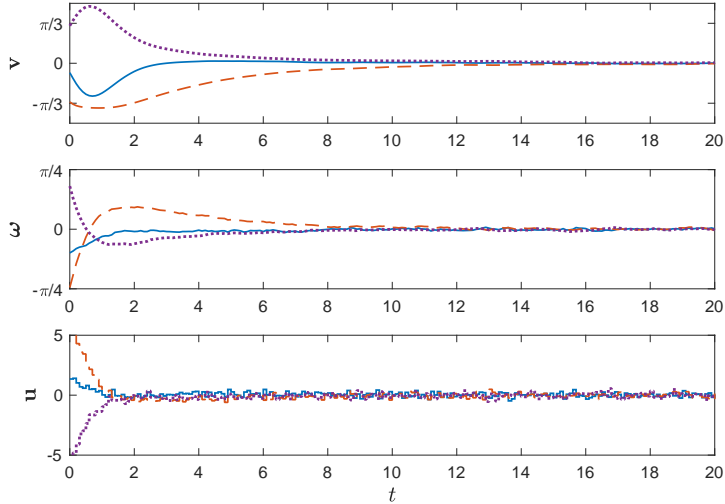}
\caption{Sample closed-loop trajectories of the rigid body system controlled by a NN feedback controller implemented with a zero-order-hold and subject to measurement noise. Solid: $\phi$, $\omega_1$, and $u_1$. Dashed: $\theta$, $\omega_2$, and $u_2$. Dotted: $\psi$, $\omega_3$, and $u_3$.}
\label{fig: satellite trajectories}
\end{figure}

In this section we perform numerical simulations of the rigid body dynamics, demonstrating that  the NN feedback controller is capable of stabilizing the system. Using \cref{eq: NN feedback control} we calculate the optimal feedback control law
\begin{equation}
\label{eq: satellite optimal control}
\bm u^{\text{NN}} (\bm v, \bm \omega)
	= - \frac{1}{W_3} \left[ \bm J^{-1} \bm B \right]^T V^{\text{NN}}_{\bm \omega} (\bm v, \bm \omega) .
\end{equation}
Recall that because we are using a time-independent value function model, the control is implemented as time-independent moving horizon feedback. Since $\bm J$ and $\bm B$ are constant matrices, we pre-compute the product $- \left[ \bm J^{-1} \bm B \right]^T / W_3$. Hence evaluation of the control requires only a forward pass through the computational graph of $V^{\text{NN}}_{\bm \omega} (\cdot)$ and a matrix multiplication.

In \cref{fig: satellite trajectories}, we plot a typical closed-loop trajectory starting from a randomly sampled initial condition. To make the simulation more realistic, we implement the controller using a zero-order-hold with a sample rate of 10 [Hz]. In addition, we corrupt inputs to the controller with Gaussian white noise with standard deviation $\sigma = 0.01 \pi$. That is, for all $t \in [t_k, t_k + 0.1]$, we apply the control
$$
\bm u (t) = \bm u^{\text{NN}} \left( \hat{\bm v} (t_k), \hat{ \bm \omega} (t_k) \right) ,
$$
where
$$
\begin{pmatrix}
\hat{\bm v} (t_k) \\ \hat{\bm \omega} (t_k)
\end{pmatrix}
\coloneqq \begin{pmatrix}
\bm v (t_k) \\
\bm \omega (t_k)
\end{pmatrix}
+ \bm n (t_k) ,
\qquad
\bm n (t_k) \sim \mathcal N \left( \bm 0, \sigma^2 \bm I \right) .
$$
In spite of this, the NN controller successfully stabilizes the system. Furthermore, the total cost of the closed-loop trajectory is $J \left[ \bm u^{\text{NN}} (\cdot) \right] = 12.67$, about $1 \%$ more than the optimal cost $J \left[ \bm u^* (\cdot) \right] = 12.52$. For comparison, a linear quadratic regulator (LQR) for \cref{eq: satellite OCP} accumulates a total cost of $J \left[ \bm u^{\text{LQR}} (\cdot) \right] = 15.95$, which is $27\%$ more than the optimal cost. Finally, short computation time is critical for implementation in real systems, and this is achieved here as each evaluation of the control takes only a couple milliseconds on both an NVIDIA RTX 2080Ti GPU and a 2012 MacBook Pro.

\section{Application to control of Burgers'-type PDE}
\label{sec: Burgers}

In this section, we test our method on high-dimensional nonlinear systems arising from a Chebyshev pseudospectral (PS) discretization of a one-dimensional forced Burgers'-type PDE. An infinite-horizon version of this problem is studied in \cite{Kalise2018}, in which the value function is approximated using a polynomial Galerkin technique. We note that in \cite{Kalise2018}, separability of the nonlinear dynamics is required to compute the high-dimensional integrals necessary in the Galerkin formulation. Our method does not require this restriction, although it does apply to this problem.

As in \cite{Kalise2018}, let $X (t, \xi) : [0, t_f] \times [-1, 1] \to \mathbb R$ satisfy the following one-dimensional controlled PDE with Dirichlet boundary conditions:
\begin{equation}
\label{eq: Burgers PDE}
\left \{ \begin{array}{ll}
X_t = X X_{\xi} + \nu X_{\xi \xi} + \alpha X e^{\beta X} + I_{\Omega} (\xi) u, & t > 0 , \xi \in (-1, 1), \\
X(t, -1) = X(t, 1) = 0 , & t > 0 , \\
X(0, \xi) = X_0 , & \xi \in (-1,1) .
\end{array} \right.
\end{equation}
For notational convenience we have written $X = X(t, \xi)$, and as before we denote $X_t = \del X / \del t$, $X_\xi = \del X / \del \xi$, and $X_{\xi \xi} = \del^2 X / \del \xi^2$. The scalar-valued control $u (t, X)$ is actuated only on $\Omega$, the support of the indicator function
$$
I_{\Omega} (\xi) \coloneqq \begin{dcases*}
1, & $\xi \in \Omega$, \\
0, & $\xi \not \in \Omega$ .
\end{dcases*}
$$
The PDE-constrained optimal control problem is
\begin{equation}
\label{eq: original Burgers OCP}
\left \{
\begin{array}{cl}
\underset{u (\cdot)}{\text{minimize}} & J \left[ u (\cdot) \right] = \displaystyle \int_t^{t_f} \mathcal L (X, u) d\tau + \frac{W_2}{2} \left \Vert X (t_f, \xi) \right \Vert_{L^2_{(-1, 1)}}^2 , \\
\text{subject to} & X_t = X X_{\xi} + \nu X_{\xi \xi} + \alpha X e^{\beta X} + I_{\Omega} (\xi) u , \\
	& X(\tau, -1) = X(\tau, 1) = 0 .
\end{array}
\right .
\end{equation}
Here
$$
\left \Vert X (\tau, \xi) \right \Vert^2_{L_{(-1,1)}^2}
\coloneqq \int_{-1}^1 X^2 (\tau, \xi) d\xi ,
\qquad
\mathcal L (X, u) = \frac{1}{2} \left \Vert X (\tau, \xi) \right \Vert_{L^2_{(-1, 1)}}^2 + \frac{W_1}{2} u^2 (\tau, X) ,
$$
and we set
$$
\begin{array}{ccccccc}
\Omega = (-0.5, -0.2), & \nu = 0.2, & \alpha = 1.5, & \beta = -0.1, & W_1 = 0.1, & W_2 = 1, & t_f = 8.
\end{array}
$$
In this problem, the goal of stabilizing $X (t, \xi)$ is made more challenging by the added reaction term, $\alpha X e^{\beta X}$, which renders the origin unstable. This can be seen clearly in \cref{fig: Burgers uncontrolled}.

To solve \cref{eq: original Burgers OCP} using our framework, we perform Chebyshev PS collocation to transform the PDE \cref{eq: Burgers PDE} into a system of ordinary differential equations (ODEs). Following \cite{Trefethen2000}, let
$$
\xi_j = \cos (j \pi / N_c),
\qquad
j = 0, 1, \dots N_c ,
$$
where $N_c+1$ is the number of collocation points. After accounting for boundary conditions, we collocate $X (t, \xi)$ at internal (non-boundary) Chebyshev points, $\xi_j$, $j = 1, 2, \dots , n$, where $n = N_c - 1$. The discretized state is defined as
$$
\bm x (t) \coloneqq \begin{pmatrix}
X (t, \xi_1), &
X (t, \xi_2), &
\dots , &
X (t, \xi_n)
\end{pmatrix}^T
: [0, t_f] \to \mathbb R^n ,
$$
and the PDE \cref{eq: Burgers PDE} becomes a system of ODEs in $n$ dimensions:
$$
\dot{\bm x} = \bm x \odot \bm {Dx} + \nu \bm D^2 \bm x + \alpha \bm x \odot e^{\beta \bm x} + \mathbb I_{\Omega} u ,
$$
In the above, ``$\odot$'' denotes element-wise multiplication (Hadamard product), $\mathbb I_{\Omega}$ is the discretized indicator function, and $\bm D, \bm D^2 \in \mathbb R^{n \times n}$ are the internal parts of the first and second order Chebyshev differentiation matrices, which are obtained by deleting the first and last rows and columns of the full matrices. This discretization automatically enforces the boundary conditions. Finally, since $X (t, \xi)$ is collocated at Chebyshev nodes, the inner product appearing in the cost function is conveniently approximated by Clenshaw-Curtis quadrature \cite{Trefethen2000}:
$$
\left \Vert X (\tau, \xi) \right \Vert^2_{L_{(-1,1)}^2}
= \int_{-1}^1 X^2 (\tau, \xi) d\xi
\approx \bm w^T \bm x^2 (\tau) ,
$$
where $\bm w \in \mathbb R^n$ are the internal Clenshaw-Curtis quadrature weights and $\bm x^2 (t)$ is calculated element-wise, i.e. $\bm x^2 = \bm x \odot \bm x$. Now the original OCP \cref{eq: original Burgers OCP} can be reformulated as an ODE-constrained problem,
\begin{equation}
\label{eq: PS Burgers OCP}
\left \{ \begin{array}{cl}
\underset{u (\cdot)}{\text{minimize}} &  \displaystyle \int_t^{t_f} \frac{1}{2} \left[ \bm w^T \bm x^2 (\tau) + W_1 u^2 (\tau, \bm x) \right] d\tau + \frac{W_2}{2} \bm w^T \bm x^2 (t_f) , \\
\text{subject to} & \dot{\bm x} = \bm x \odot \bm {Dx} + \nu \bm D^2 \bm x + \alpha \bm x \odot e^{\beta \bm x} + \mathbb I_{\Omega} u .
\end{array} \right.
\end{equation}

\subsection{Learning high-dimensional value functions}
\label{sec: Burgers results}

The state dimension $n$ of the OCP \cref{eq: PS Burgers OCP} can be adjusted, presenting a good opportunity to test the scalability of our algorithms. For this problem, we learn the value function $V = V (t, \bm x)$ with time-dependence, rather than just $V (0, \bm x)$ as in \cref{sec: rigid body}. Consequently, the resulting controls can be implemented as time-dependent controls or with a moving horizon. We consider the following domain of initial conditions:
\begin{equation}
\label{eq: burgers initial conditions}
\mathbb X_0 = \left \{ \left. \bm x \in \mathbb R^n \right| -2 \leq x_j \leq 2, j = 1, 2, \dots, n \right \} .
\end{equation}

Using the proposed adaptive deep learning framework, we approximate solutions to \cref{eq: PS Burgers OCP} in $n=10$, 20, and 30 dimensions. We focus on demonstrating what is possible using our approach, rather than carrying out a detailed study of its effectiveness under different parameter tunings. In \cite{Kalise2018} an infinite-horizon version of the problem is solved up to twelve dimensions, but the accuracy of the solution is not readily verifiable. The ability to conveniently measure model accuracy for general high-dimensional problems with \textit{no known analytical solution} is a key advantage of our framework.

For each discretized OCP, $n=10$, 20, and 30, we apply the time-marching strategy to build an initial training data set $\mathcal D_{\text{train}}^1$ from 30 uniformly sampled initial conditions, $\bm x_0^{(i)} \in \mathbb X_0$, $i = 1, 2, \dots, 30$. For each initial condition $\bm x_0^{(i)}$,  the BVP solver outputs an optimal trajectory $\left \{ \bm x^{(i)} \left( t_k \right) \right \}$, evaluated at collocation points $t_k \in [0, t_f]$ chosen by the solver. Typically this can be a few hundred per initial condition, depending on the state dimension $n$ and the BVP solver tolerances. Since these data sets can be get quite large, we often train on randomly selected subsets of the data. This can significantly improve training speed without sacrificing accuracy. When neeeded, we solve additional BVPs to expand the data set as described in \cref{sec: adaptive sampling implementation}. We use the same NN architecture as in \cref{sec: rigid body}, with three hidden layers with 64 neurons each. We set $C = 0.3$, 1.3, and 1.8 for $n = 10$, 20, and 30, respectively, and use $\mu = 10$ in all cases.

\begin{table}[t!]
\centering
\begin{tabular}{| c | c | c | c | c |}
\hline
&&&&\\[-1em]
$n$ & num. trajectories & training time & value RMAE & gradient RM$L^2$ \\
\hline
&&&&\\[-1em]
10 & 163 & 25 min & $1.3 \times 10^{-3}$ & $5.7 \times 10^{-3}$ \\
20 & 128 & 48 min & $5.0 \times 10^{-3}$ & $1.1 \times 10^{-2}$ \\
30 & 145 & 62 min & $1.3 \times 10^{-2}$ & $2.2 \times 10^{-2}$ \\
\hline
\end{tabular}
\caption{Validation accuracy of NNs for solving the collocated Burgers'-type OCP \cref{eq: PS Burgers OCP}, depending on the state dimension $n$. Training time includes time spent generating additional data according to \cref{alg: training}. All NNs are trained on an NVIDIA RTX 2080Ti GPU.}
\label{table: Burgers accuracy results}
\end{table}

In \cref{table: Burgers accuracy results}, we present validation accuracy results for the trained NNs. We include the RMAE in predicting the value function and the RM$L^2$ error in predicting the costate, $\bm \lambda (t; \bm x_0) \approx V^{\text{NN}}_{\bm x} \left( t, \bm x (t; \bm x_0) \right)$. Accuracy is measured empirically on independently generated validation data sets comprised of trajectories from 50 randomly selected initial conditions. We find that the trained NNs accurately predict both the value function and its gradient, even in 30 dimensions.

\Cref{table: Burgers accuracy results} also shows the total number of sample trajectories seen by the NN, including the initial data $\mathcal D_{\text{train}}^1$. It may seem surprising that we are able to reach the same level of accuracy in higher dimensions with similar numbers of sample trajectories. This happens because the BVP solver usually needs more collocation points for larger problems, thus producing more data per trajectory. Consequently, fewer trajectories are needed to fulfill the data set size recommendation \cref{eq: our sample size selection}. Similarly, in \cref{sec: rigid body} we used data only for $t=0$, so we needed thousands of trajectories to fill the state space. This suggests that learning the time-dependent value function can be more efficient than learning $V(0, \bm x)$ only. Note that, if preferred, the time-dependent controller can still be used with a moving horizon like in \cref{sec: rigid body}.

Lastly, \cref{table: Burgers accuracy results} shows the training time for each NN, including time spent testing convergence and generating additional trajectories on the fly, but not time spent generating the initial data. Generating data becomes the most expensive computation as $n$ increases, but even so we find that computational effort scales reasonably with the problem dimension. Furthermore, it is possible to obtain a rough low-fidelity NN model in just minutes as shown in \cref{table: Burgers BVP convergence NN}, which in turn allows for more efficient data generation. This demonstrates the viability of the proposed method for solving high-dimensional optimal control problems.

\subsection{NN warm start for fast and reliable BVP solutions}
\label{sec: Burgers adaptive sampling}

\begin{table}[t!]
\centering
\begin{tabular}{| c | c | c | c |}
\hline
&&&\\[-1em]
$n$ & $K$ & $\%$ BVP convergence & mean integration time \\
\hline
&&&\\[-1em]
\multirow{3}{*}{10} & 4 & $40 \%$ & 0.7 s \\
& 6 & $83 \%$ & 0.8 s \\
& 10 & $90 \%$ & 1.3 s \\
\hline
&&&\\[-1em]
\multirow{3}{*}{20} & 4 & $46 \%$ & 3.6 s \\
& 5 & $86 \%$ & 4.2 s \\
& 6 & $99 \%$ & 4.7 s \\
\hline
&&&\\[-1em]
\multirow{3}{*}{30} & 4 & $47 \%$ & 11.3 s \\
& 6 & $90 \%$ & 14.6 s \\
& 8 & $100 \%$ & 19.1 s \\
\hline
\end{tabular}
\caption{Convergence of BVP solutions for \cref{eq: PS Burgers OCP} when using the time-marching trick, depending on the problem dimension, $n$, and the number of steps in the sequence $\{ t_k \}_{k=1}^K$. BVP integration time is measured only on successful attempts.}
\label{table: Burgers BVP convergence}
\begin{tabular}{| c | c | c | c | c | c |}
\hline
&&&&&\\[-1em]
$n$ & $\mu$ & training time & gradient RM$L^2$ & $\%$ BVP conv. & mean int. time \\
\hline
&&&&&\\[-1em]
\multirow{3}{*}{10} & $10^{-8}$ & 20 s & $2.1 \times 10^{-1}$ & $96 \% $ & 0.8 s \\
& $10^{-4}$ & 31 s & $9.8 \times 10^{-2}$ & $99 \% $ & 0.8 s \\
& 1 & 56 s & $4.2 \times 10^{-2}$ & $88 \%$ & 0.6 s \\
\hline
&&&&&\\[-1em]
\multirow{3}{*}{20} & $10^{-8}$ & 29 s & $3.7 \times 10^{-1}$ & $74 \% $ & 2.9 s \\
& $10^{-4}$ & 47 s & $9.8 \times 10^{-2}$ & $91 \% $ & 2.5 s \\
& 1 & 76 s & $6.5 \times 10^{-2}$ & $98 \%$ & 2.5 s \\
\hline
&&&&&\\[-1em]
\multirow{3}{*}{30} & $10^{-8}$ & 38 s & $3.0 \times 10^{-1}$ & $79 \%$ & 7.1 s \\
& $10^{-4}$ & 125 s & $7.6 \times 10^{-2}$ & $94 \% $ & 6.9 s \\
& 1 & 189 s & $7.4 \times 10^{-2}$ & $96 \%$ & 7.1 s \\
\hline
\end{tabular}
\caption{Convergence of BVP solutions for \cref{eq: PS Burgers OCP} when using NN warm start with NNs of varying gradient prediction accuracy. BVP integration time is measured only on successful attempts.}
\label{table: Burgers BVP convergence NN}
\end{table}

In our experience, generating the initial training data set can be the most computationally demanding part of the process, especially as the problem dimension $n$ increases. Consequently, for difficult high-dimensional problems it may be impractical to generate a large-enough data set from scratch. This obstacle can be largely overcome by using partially-trained/low-fidelity NNs to aid in further data generation. In this section, we briefly compare the reliability and speed of BVP convergence between our two strategies: time-marching and NN warm start. These experiments demonstrate the importance of NN guesses for high-dimensional data generation.

For each of $n = 10$, 20, and 30, we randomly sample a set of 1000 candidate points from the domain $\mathbb X_0$ defined in \cref{eq: burgers initial conditions}. From these we choose 100 points with the largest predicted value gradient. The set of initial conditions is fixed for each $n$. Next we proceed as in \cref{sec: satellite adaptive sampling}, solving each BVP by time-marching with different number of time-marching iterations, $K$. We tune each time sequence to improve convergence as much as possible. Results are summarized in \cref{table: Burgers BVP convergence}. We then solve the same BVPs directly over the whole time interval $t \in [0, 8]$ with NN warm start. These NNs are trained on fixed data sets containing only 30 trajectories, but with different gradient loss weights $\mu$, resulting in varying costate prediction accuracy. We also limit the number of L-BFGS iterations so that each model is trained only for a short time. Results are given in \cref{table: Burgers BVP convergence NN}.

As before, we find that even NNs with relatively large costate prediction error enable consistently convergent BVP solutions. Time-marching also works once the sequence of time steps $\{ t_k \}_{k=1}^K$ is properly tuned, but the speed of this method scales poorly with $n$. Now the advantage of utilizing NNs to aid in data generation becomes clear: when $n$ is large, the average time needed for convergence when using NN warm start is drastically lower than that of the time-marching trick. This approach also requires no tuning of the time-marching sequence. Because low-fidelity NNs are quick to train, training such a NN and then using it to aid in data generation is the most efficient strategy for building larger data sets.

\subsection{Closed-loop simulations}

\begin{figure}[t!]
\centering

\begin{subfigure}{\textwidth}
\centering
\includegraphics[width = .49\textwidth]{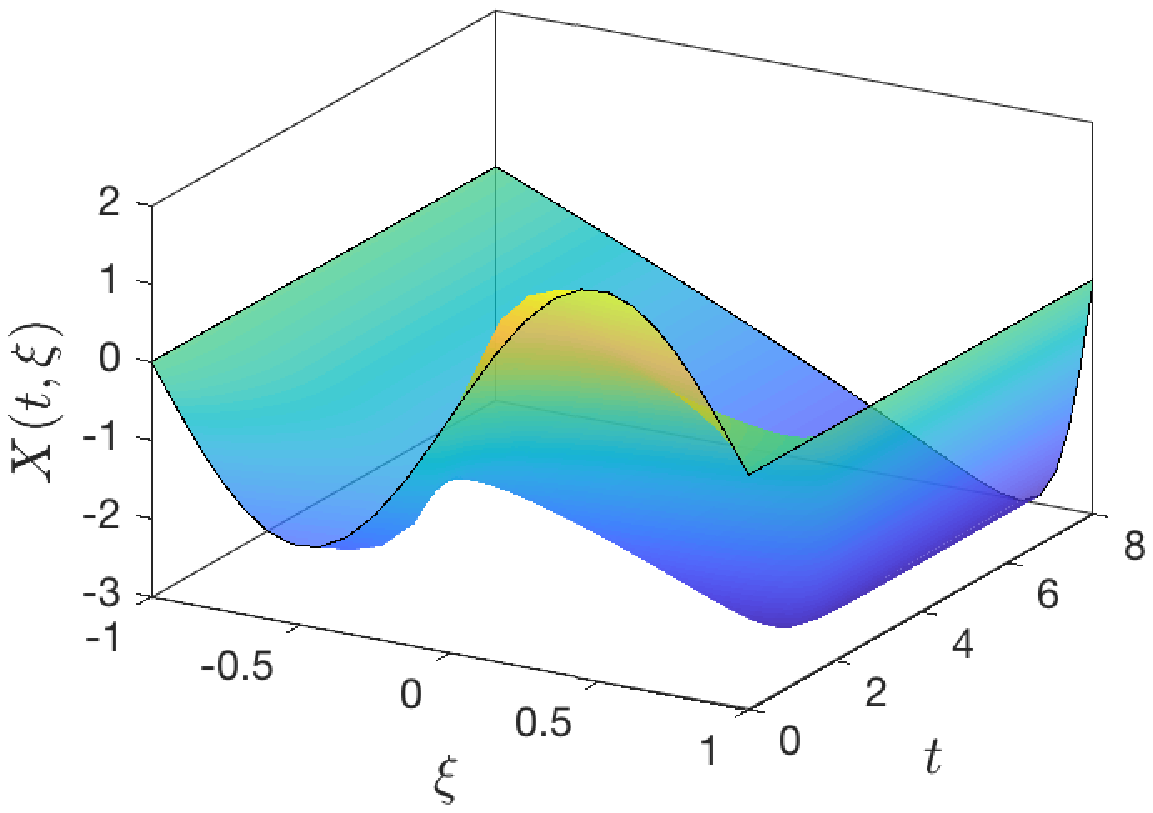}
\includegraphics[width = .49\textwidth]{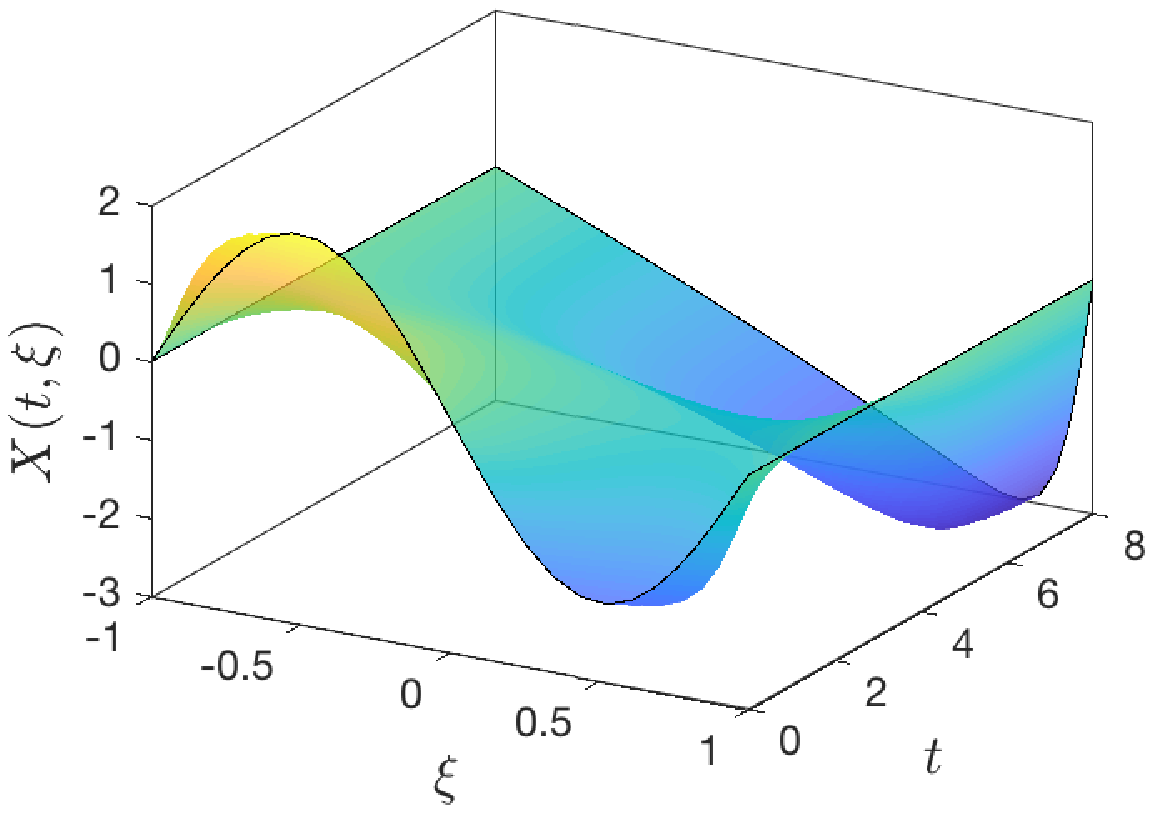}
\caption{Uncontrolled dynamics.}
\label{fig: Burgers uncontrolled}
\end{subfigure}

\begin{subfigure}{\textwidth}
\centering
\includegraphics[width = .49\textwidth]{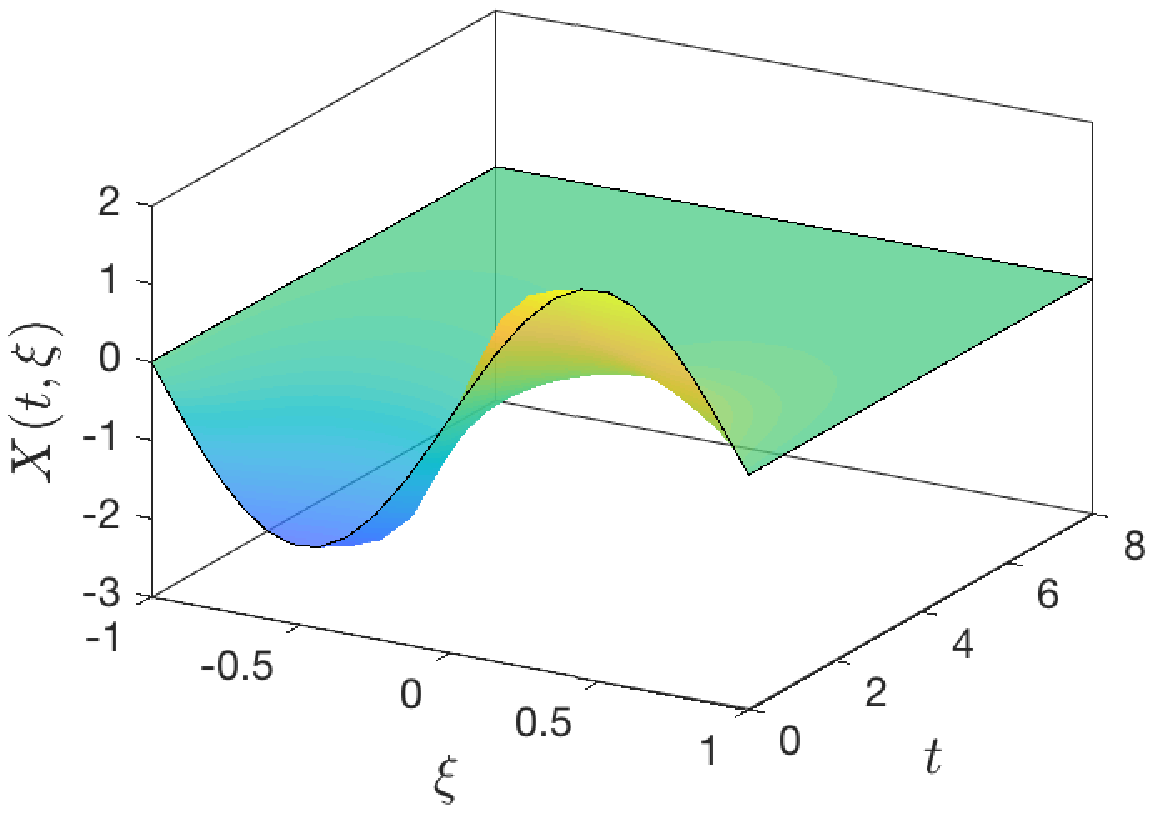}
\includegraphics[width = .49\textwidth]{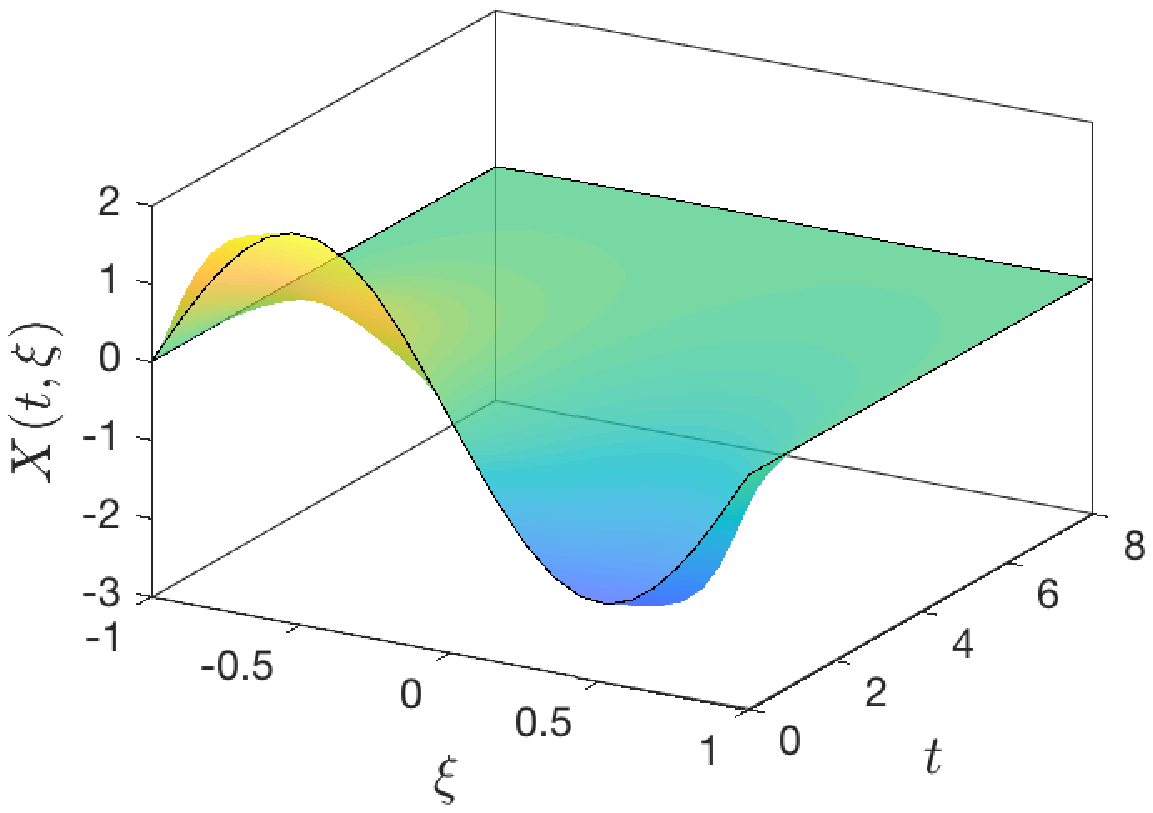}
\caption{NN-controlled dynamics.}
\label{fig: Burgers controlled}
\end{subfigure}

\begin{subfigure}{\textwidth}
\centering
\includegraphics[width = .49\textwidth]{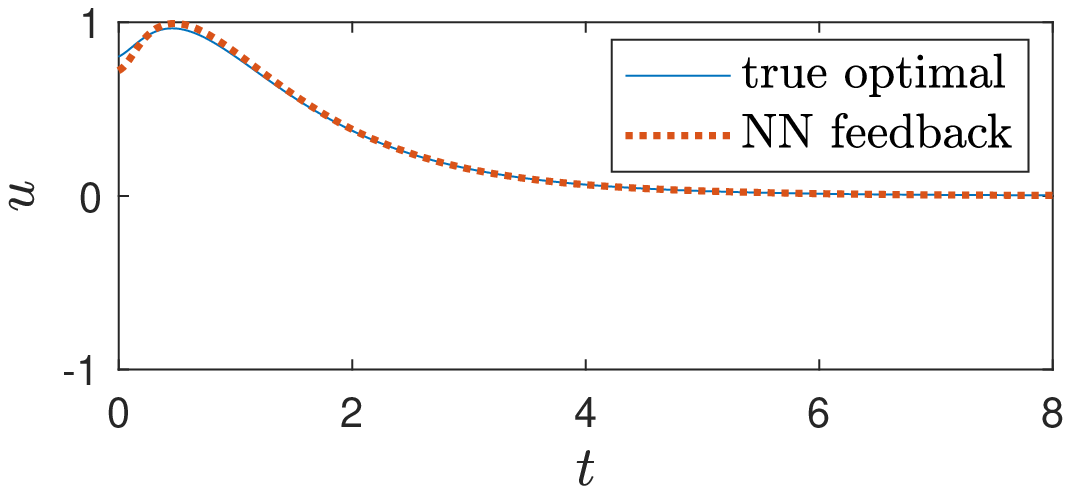}
\includegraphics[width = .49\textwidth]{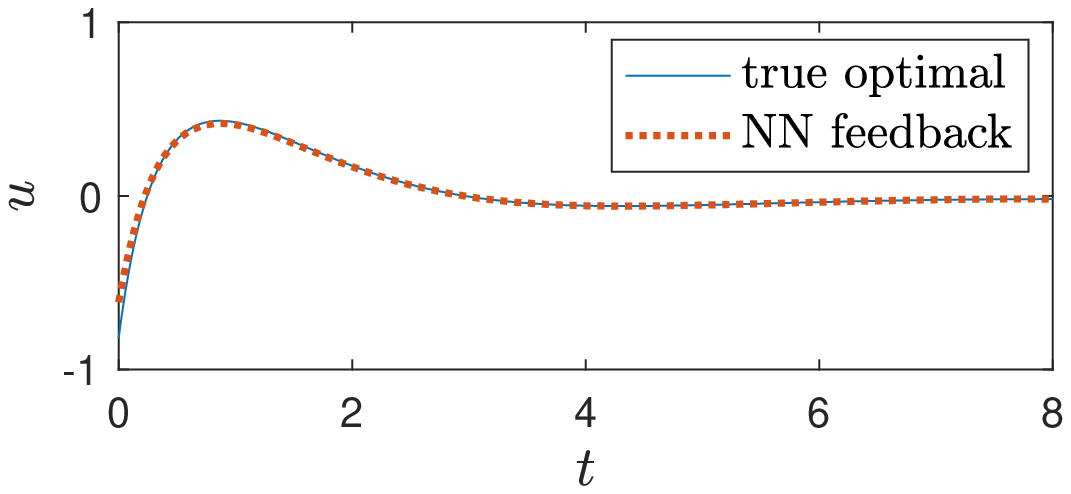}
\caption{Comparison of true optimal control (open-loop BVP solution) and NN control profiles.}
\label{fig: Burgers control}
\end{subfigure}

\caption{Simulations of the collocated Burgers'-type PDE \cref{eq: Burgers PDE} in $n=30$ dimensions. Left column: $X(0,\xi) = 2 \sin (\pi \xi)$. Right column: $X(0,\xi) = -2 \sin (\pi \xi)$.}
\label{fig: Burgers dynamics}
\end{figure}

In this section we use simulations to demonstrate that the feedback control output by the trained NN not only stabilizes the high-dimensional system, but that it is close to the true optimal control. The optimal feedback control law can again be calculated with \cref{eq: NN feedback control}, from which we obtain
\begin{equation}
\label{eq: Burgers optimal control}
u^{\text{NN}} (t, \bm x) = - \frac{1}{W_1} \mathbb I_{\Omega}^T \, V_{\bm x}^{\text{NN}} (t, \bm x) .
\end{equation}

In \cref{fig: Burgers dynamics}, we plot the uncontrolled (\cref{fig: Burgers uncontrolled}) and closed-loop controlled dynamics (\cref{fig: Burgers controlled}), starting from two different initial conditions, $X(0,\xi) = 2\sin(\pi \xi)$ and $X(0,\xi) = -2\sin(\pi \xi)$, where the dimension of the discretized system is $n=30$. For both of these initial conditions (and almost all others tested), the NN controller successfully stabilizes the open-loop unstable origin. Further, as shown in \cref{fig: Burgers control}, the NN-generated controls are very close to the true optimal controls which are calculated by solving the associated BVPs. Finally, the speed of online control computation is not sensitive to the problem dimension: each evaluation still takes just milliseconds on both an NVIDIA RTX 2080Ti GPU and a 2012 MacBook Pro.

\section{Conclusion}

In this paper, we have developed a novel machine learning framework for solving HJB equations and designing candidate optimal feedback controllers. Unlike many other state of the art techniques, our method does not require finite difference approximations of the gradient nor strict restrictions on the structure of the dynamics. The causality-free algorithm we use for data generation enables application to high-dimensional systems and validation of model accuracy. We also emphasize that while our method is data-driven, by leveraging the costate data we are able to train more physically-consistent models and better controllers with surprisingly small data sets.

The proposed method is not only a consumer of data, but through adaptive data generation it can also be used build rich data sets with points anywhere in a semi-global domain. Thus the value function and control are valid for large ranges of dynamic states, rather than just in the neighborhood of some nominal trajectory. Furthermore, data can be generated specifically near complicated regions of the value function or where the required control effort is large. This in turn allows us to train more accurate NN models and synthesize better-performing controllers, or employ other data-driven methods.

We have demonstrated the possibility for use of the framework in a practical setting by designing candidate optimal feedback controllers for a six-dimensional nonlinear rigid body. The potential for scalability of the method is demonstrated by solving HJB equations in up to 30 dimensions using limited data, and empirical validation indicates that the NN models are good approximations of the value function. How well the proposed techniques work for even larger problems remains an open question. Indeed, understanding the scalability of deep learning methods in general is still an active area of research. Nevertheless, we are encouraged by the simulations in \cref{sec: Burgers} which suggest that the method may scale quite well for moderately high-dimensional problems. In addition, the computational burden associated with an increase in dimensionality is incurred entirely offline: due to the structure of NNs, increasing the dimension has a negligible impact on the speed of online control calculation.

These promising results leave plenty of room for future development. Of special interest are extensions of the framework to solve problems with free final time and state and control constraints, which appear ubiquitously in practical applications. Such problems typically give rise to non-unique solutions of PMP and non-smooth value functions, thus presenting substantial challenges for both data generation and neural network modeling. Overcoming these obstacles would open the door to solving many important and difficult optimal control problems.

\bibliography{bibliography_HJB_SIAM}

\begin{thebibliography}{10}

\bibitem{TensorFlow}
{\sc M.~Abadi, A.~Agarwal, P.~Barham, et~al.}, {\em {TensorFlow}: Large-scale
  machine learning on heterogeneous systems}, 2016,
  \url{https://arxiv.org/abs/1603.04467}.

\bibitem{Khalaf2005}
{\sc M.~Abu-Khalaf and F.~L. Lewis}, {\em Nearly optimal control laws for
  nonlinear systems with saturating actuators using a neural network {HJB}
  approach}, Automatica, 41 (2005), pp.~779--791,
  \url{https://doi.org/10.1016/j.automatica.2004.11.034}.

\bibitem{Albrekht1961}
{\sc E.~Al'brekht}, {\em On the optimal stabilization of nonlinear systems}, J.
  Appl. Math. Mech., 25(5) (1961), pp.~1254--1266,
  \url{https://doi.org/10.1016/0021-8928(61)90005-3}.

\bibitem{Bachouch2018}
{\sc A.~Bachouch, C.~Hur\'e, N.~Langren\'e, and H.~Pham}, {\em Deep neural
  networks algorithms for stochastic control problems on finite horizon:
  numerical applications}, 2018, \url{https://arxiv.org/abs/1812.05916}.

\bibitem{Bokanowski2013}
{\sc O.~Bokanowski, J.~Garcke, M.~Griebel, and I.~Klompmaker}, {\em An adaptive
  sparse grid semi-{L}agrangian scheme for first order {Hamilton-Jacobi
  Bellman} equations}, J. Sci. Comput., 55 (2013), pp.~575--605,
  \url{https://doi.org/10.1007/s10915-012-9648-x}.

\bibitem{Bottou2018}
{\sc L.~Bottou, F.~E. Curtis, and J.~Nocedal}, {\em Optimization methods for
  large-scale machine learning}, SIAM Rev., 60 (2018), pp.~223--311,
  \url{https://doi.org/10.1137/16M1080173}.

\bibitem{Byrd2012}
{\sc R.~H. Byrd, G.~M. Chin, J.~Nocedal, and Y.~Wu}, {\em Sample size selection
  in optimization methods for machine learning}, Math. Program., 134 (2012),
  pp.~127--155, \url{https://doi.org/10.1007/s10107-012-0572-5}.

\bibitem{Byrd1995}
{\sc R.~H. Byrd, P.~Lu, J.~Nocedal, and C.~Zhu}, {\em A limited memory
  algorithm for bound constrained optimization}, SIAM J. Sci. Comput., 16
  (1995), pp.~1190--1208, \url{https://doi.org/10.1137/0916069}.

\bibitem{Cacace2012}
{\sc S.~Cacace, E.~Cristiani, M.~Falcone, and A.~Picarelli}, {\em A patchy
  dynamic programming scheme for a class of {H}amilton-{J}acobi-{B}ellman
  equations}, SIAM J. Sci. Comput., 34 (2012), pp.~A2625--A2649,
  \url{https://doi.org/10.1137/110841576}.

\bibitem{Cheng2007}
{\sc T.~Cheng, F.~L. Lewis, and M.~Abu-Khalaf}, {\em
  Fixed-final-time-constrained optimal control of nonlinear systems using
  neural network {HJB} approach}, {IEEE} Trans. Neural Netw., 18 (2007),
  pp.~1725--1737, \url{https://doi.org/10.1109/TNN.2007.905848}.

\bibitem{Chow2019}
{\sc Y.~T. Chow, J.~Darbon, S.~Osher, and W.~Yin}, {\em Algorithm for
  overcoming the curse of dimensionality for state-dependent {Hamilton-Jacobi}
  equations}, J. Comput. Phys., 387 (2019), pp.~376--409,
  \url{https://doi.org/10.1016/j.jcp.2019.01.051}.

\bibitem{Crandall1983}
{\sc M.~G. Crandall and P.-L. Lions}, {\em Viscosity solutions of
  {H}amilton-{J}acobi equations}, Trans. Amer. Math. Soc., 277 (1983),
  pp.~1--42, \url{https://doi.org/10.2307/1999343}.

\bibitem{Darbon2020}
{\sc J.~Darbon, G.~P. Langlois, and T.~Meng}, {\em Overcoming the curse of
  dimensionality for some {H}amilton--{J}acobi partial differential equations
  via neural network architectures}, Res. Math. Sci., 7 (2020), p.~20,
  \url{https://doi.org/10.1007/s40687-020-00215-6}.

\bibitem{Darbon2021}
{\sc J.~Darbon and T.~Meng}, {\em On some neural network architectures that can
  represent viscosity solutions of certain high dimensional
  {H}amilton--{J}acobi partial differential equations}, J. Comput. Phys., 425
  (2021), p.~109907,
  \url{https://doi.org/https://doi.org/10.1016/j.jcp.2020.109907}.

\bibitem{Darbon2016}
{\sc J.~Darbon and S.~Osher}, {\em Algorithms for overcoming the curse of
  dimensionality for certain {H}amilton-{J}acobi equations arising in control
  theory and elsewhere}, Res. Math. Sci., 3 (2016),
  \url{https://doi.org/10.1186/s40687-016-0068-7}.

\bibitem{Diebel2006}
{\sc J.~Diebel}, {\em Representing attitude: {E}uler angles, unit quaternions,
  and rotation vectors}, 2006,
  \url{{https://www.astro.rug.nl/software/kapteyn-beta/\_downloads/attitude.pdf}}
  (accessed 2020-05-16).

\bibitem{Falcone2013}
{\sc M.~Falcone and R.~Ferretti}, {\em Semi-{Lagrangian} Approximation Schemes
  for Linear and {Hamilton-Jacobi} Equations}, Society for Industrial and
  Applied Mathematics, Philadelphia, PA, 2013,
  \url{https://doi.org/10.1137/1.9781611973051}.

\bibitem{Han2018_PNAS}
{\sc J.~Han, A.~Jentzen, and W.~E}, {\em Solving high-dimensional partial
  differential equations using deep learning}, Proc. Natl. Acad. Sci. USA, 115
  (2018), pp.~8505--8510, \url{https://doi.org/10.1073/pnas.1718942115}.

\bibitem{Hure2018}
{\sc C.~Hur\'e, H.~Pham, A.~Bachouch, and N.~Langren\'e}, {\em Deep neural
  networks algorithms for stochastic control problems on finite horizon, part
  {I}: convergence analysis}, 2018, \url{https://arxiv.org/abs/1812.04300}.

\bibitem{Izzo2019}
{\sc D.~Izzo, E.~\"{O}zt\"{u}rk, and M.~M\"{a}rtens}, {\em Interplanetary
  transfers via deep representations of the optimal policy and/or of the value
  function}, in Genetic and Evolutionary Computation Conference, 2019,
  pp.~1971---1979, \url{https://doi.org/10.1145/3319619.3326834}.

\bibitem{Jiang2016}
{\sc F.~Jiang, G.~Chou, M.~Chen, and C.~J. Tomlin}, {\em Using neural networks
  to compute approximate and guaranteed feasible {Hamilton-Jacobi-Bellman PDE}
  solutions}, 2016, \url{https://arxiv.org/abs/1611.03158}.

\bibitem{Kalise2018}
{\sc D.~Kalise and K.~Kunisch}, {\em Polynomial approximation of
  high-dimensional {Hamilton-Jacobi-Bellman} equations and applications to
  feedback control of semilinear parabolic {PDE}s}, SIAM J. Sci. Comput., 40
  (2018), pp.~A629--A652, \url{https://doi.org/10.1137/17M1116635}.

\bibitem{Kang1992}
{\sc W.~Kang, P.~De, and A.~Isidori}, {\em Flight control in a windshear via
  nonlinear $h_\infty$ methods}, in Proceedings of the 31st IEEE Conference on
  Decision and Control, vol.~1, 1992, pp.~1135--1142.

\bibitem{Kang2015}
{\sc W.~Kang and L.~C. Wilcox}, {\em A causality free computational method for
  {HJB} equations with application to rigid body satellites}, in AIAA Guidance,
  Navigations, and Control Conference, 2015, pp.~1--10,
  \url{https://doi.org/10.2514/6.2015-2009}.

\bibitem{Kang2017_COA}
{\sc W.~Kang and L.~C. Wilcox}, {\em Mitigating the curse of dimensionality:
  Sparse grid characteristics method for optimal feedback control and {HJB}
  equations}, Comput. Optim. Appl., 68 (2017), pp.~289--315,
  \url{https://doi.org/10.1007/s10589-017-9910-0}.

\bibitem{Kierzenka2001}
{\sc J.~Kierzenka and L.~F. Shampine}, {\em A {BVP} solver based on residual
  control and the {MATLAB PSE}}, ACM Trans. Math. Softw., 27 (2001),
  pp.~299--316, \url{https://doi.org/10.1145/502800.502801}.

\bibitem{Liberzon2011}
{\sc D.~Liberzon}, {\em Calculus of Variations and Optimal Control Theory: A
  Concise Introduction}, Princeton University Press, Princeton, NJ, 2011,
  \url{https://doi.org/10.2307/j.ctvcm4g0s}.

\bibitem{Lukes1969}
{\sc D.~Lukes}, {\em Optimal regulation of nonlinear dynamical systems}, SIAM
  J. Control, 7 (1969), pp.~75--100, \url{https://doi.org/10.1137/0307007}.

\bibitem{Mangasarian1966}
{\sc O.~L. Mangasarian}, {\em Sufficient conditions for the optimal control of
  nonlinear systems}, SIAM J. Control, 4 (1966), pp.~139--152,
  \url{https://doi.org/10.1137/0304013}.

\bibitem{Nakamura2020}
{\sc T.~Nakamura-Zimmerer, Q.~Gong, and W.~Kang}, {\em A causality-free neural
  network method for high-dimensional {H}amilton-{J}acobi-{B}ellman equations},
  in American Control Conference (ACC), 2020, pp.~787--793,
  \url{https://doi.org/10.23919/ACC45564.2020.9147270}.

\bibitem{Navasca2007}
{\sc C.~Navasca and A.~J. Krener}, {\em Patchy Solutions of
  {H}amilton-{J}acobi-{B}ellman Partial Differential Equations}, Springer,
  Berlin-Heidelberg, 2007, pp.~251--270,
  \url{https://doi.org/10.1007/978-3-540-73570-0_20}.

\bibitem{Osher1988}
{\sc S.~Osher and J.~A. Sethian}, {\em Fronts propagating with
  curvature-dependent speed: Algorithms based on {H}amilton-{J}acobi
  formulations}, J. Comput. Phys., 79 (1988), pp.~12--49,
  \url{https://doi.org/10.1016/0021-9991(88)90002-2}.

\bibitem{Raissi2019}
{\sc M.~Raissi, P.~Perdikaris, and G.~Karniadakis}, {\em Physics-informed
  neural networks: A deep learning framework for solving forward and inverse
  problems involving nonlinear partial differential equations}, J. Comput.
  Phys., 378 (2019), pp.~686--707,
  \url{https://doi.org/10.1016/j.jcp.2018.10.045}.

\bibitem{Sirignano2018}
{\sc J.~Sirignano and K.~Spiliopoulos}, {\em {DGM}: A deep learning algorithm
  for solving partial differential equations}, J. Comput. Phys., 375 (2018),
  pp.~1339--1364, \url{https://doi.org/10.1016/j.jcp.2018.08.029}.

\bibitem{Tassa2007}
{\sc Y.~Tassa and T.~Erez}, {\em Least squares solutions of the {HJB} equation
  with neural network value-function approximators}, {IEEE} Trans. Neural
  Netw., 18 (2007), pp.~1031--1041,
  \url{https://doi.org/10.1109/TNN.2007.899249}.

\bibitem{Trefethen2000}
{\sc L.~N. Trefethen}, {\em Spectral Methods in MATLAB}, Society for Industrial
  and Applied Mathematics, Philadelphia, PA, 2000,
  \url{https://doi.org/10.1137/1.9780898719598}.

\bibitem{SciPy}
{\sc P.~Virtanen, R.~Gommers, T.~E. Oliphant, and et. al.}, {\em {SciPy} 1.0:
  Fundamental algorithms for scientific computing in {Python}}, Nat. Methods,
  17 (2020), pp.~261--272, \url{https://doi.org/10.1038/s41592-019-0686-2}.

\bibitem{Yegorov2018}
{\sc I.~Yegorov and P.~M. Dower}, {\em Perspectives on characteristics based
  curse-of-dimensionality-free numerical approaches for solving
  {Hamilton-Jacobi} equations}, Appl. Math. Optim.,  (2018),
  \url{https://doi.org/10.1007/s00245-018-9509-6}.

\end{thebibliography}

\end{document}